\documentclass[12pt]{article}

\usepackage{latexsym}
\usepackage{amssymb}

\font\tenmsb=msbm10
\font\sevenmsb=msbm7
\font\fivemsb=msbm5

\newfam\msbfam
\textfont\msbfam=\tenmsb
\scriptfont\msbfam=\sevenmsb
\scriptscriptfont\msbfam=\fivemsb
\def\Bbb#1{{\fam\msbfam #1}}

\makeatletter
\ifnum\@ptsize=0  \fi \ifnum\@ptsize=2  \fi 

\newcommand\sE{{\cal E}}

\newcommand\sI{{\cal I}}

\newcommand\sL{{\cal L}}
\newcommand\sO{{\cal O}}

\newcommand\sC{{\cal C}}

\newcommand\zed{{\Bbb Z}}

\newcommand\bQ{{\Bbb Q}}
\newcommand\bN{{\Bbb N}}
\newcommand\bP{{\mathbb P}}

\newcommand\la{{\longrightarrow}}

\newcommand\proof{{\noindent\bf Proof.\ }} 
\newcommand\qed{{\hspace*{\fill}Q.E.D.\vskip12pt plus 1pt}}

\newtheorem{theorem}{Theorem}[section]
\newtheorem{lemma}[theorem]{Lemma}
\newtheorem{corollary}[theorem]{Corollary} 
\newtheorem{proposition}[theorem]{Proposition} 
 
\newtheorem{re}[theorem]{Remark} 
\newtheorem{setup}[theorem]{Setup}
 
\newtheorem{definition}[theorem]{Definition}

\newtheorem{example}[theorem]{Example}
\newtheorem{defprop}[theorem]{Definition-Proposition}

\newenvironment{remark}{\begin{re}\em}{\end{re}}

\begin{document}
\title {Nef reduction and anticanonical bundles} \author{Thomas Bauer and Thomas Peternell} \date{ }
\maketitle

\vspace*{-0.5in}\section*{Introduction} 

Projective manifolds $X$ with nef anticanonical bundles (i.e. $-K_X \cdot C = \det T_X \cdot C \geq 0$ for all curves 
$C \subset X$) can be regarded as an interpolation between 
Fano manifolds (ample anticanonical bundle) and Calabi-Yau manifolds resp.\ tori and symplectic manifolds
(trivial canonical bundle). A differential-geometric analogue are varieties with semi-positive
Ricci curvature although this class is strictly smaller -- to get the correct picture one has to consider
sequences of metrics and make the negative part smaller and smaller. However we will work completely in the context of algebraic
geometry.  \\
Our aim is twofold: classification and, as a consequence, boundedness in case of dimension 3. We shall not consider
threefods with trivial canonical bundles, the eventual boundedness of Calabi-Yau threefolds still being unknown.
Fano threefolds have been classified a long time ago and threefolds with big and nef anticanonical bundle are very much
related with $\bQ$-Fano threefolds; therefore we will concentrate here on {\it projective threefolds $X$ with 
$-K_X$ nef and $K_X^3 = 0,$ but $K_X \not \equiv 0.$} \\
The essential problem is to distinguish the positive and flat directions in $X.$ There are three main tools to do that:
\begin{itemize}
\item the Albanese map
\item Mori theory
\item the nef reduction.
\end{itemize}

Given a normal projective variety $X$ and a nef line bundle $L$, the nef reduction produces an almost holomorphic 
dominant
meromorphic map $f: X \rightharpoonup B$ with connected fibers such that 
\begin{enumerate}
\item $L$ is numerically trivial on all compact fibers $F$ of $f$ of dimension $\dim X - \dim B$
\item for a general point $x \in X$ and every irreducible curve $C$ passing through $x$ such that $\dim f(C) > 0,$ we have
$L \cdot C > 0.$ 
\end{enumerate} 

The number $\dim B$ is an invariant of $L$ (actually $f$ is birationally determined) and is called the nef dimension $n(L)$.
We will apply this to $L = -K_X$ and find that in case $1 \leq n(-K_X) \leq 2,$ some multiple of $-K_X$ is spanned and provides
the nef reduction (Theorem 2.1). This theorem is actually first  established in the case when $X$ is rationally connected while the other 
cases are done a posteriori with further knowledge of the structure of the variety. 
\\
The first thing in classification is of course the study of the Albanese. It is known that $\alpha: X \to {\rm Alb}(X)$ is 
a surjective submersion. Since we are interested only in classification up to finite \'etale cover, we will assume that
the irregularity $q(X) $ is maximal (with respect to finite \'etale covers) and then, possibly after another cover, 
Theorem 4.2 provides a precise structure if $n(-K_X) = 1$ or $2$ -- essentially $X$ is a product -- and allows to show that threefolds with positive irregularity are
bounded up to finite \'etale cover, also if $n(-K_X) = 3.$ \\
The Albanese theory being settled, we may now assume that $q = 0,$ even after finite \'etale cover. This means that $X$ is
simply connected. If now $X$ is not rationally connected, we can use the rational quotient and it turns out that after finite 
\'etale cover, $X$ is a product $\bP_1 \times $K3 (Theorem 3.1). \\
So we are reduced to rationally connected threefolds. Combining the holomorphic nef reduction and Mori theory which are
so to speak ``perpendicular'', we arrive at several structure theorems (sections 5 and 6) if $n(-K_X) = 1$ or $2.$ \\
The case $n(-K_X) = 3$ is studied in sect.\ 7. This condition means that $-K_X$ is ample on all irreducible members of
covering families of curves. Recalling our general assumption that $K_X^3 = 0,$ we show that $K_X^2 \ne 0$ and -- although it
has a non-zero fixed part -- the anticanonical system $\vert -K_X \vert $ induces a fibration $f: X \to \bP_1.$ 
The general fiber $F$ has $n(-K_F) = 2$ and therefore it is either $\bP_2$ blown up in 9 points but without elliptic fibration 
or a special $\bP_1$-bundle over an elliptic curve. We study in detail the first case under the genericity assumption 
that the unique element in $\vert -K_F \vert$ is a smooth elliptic curve (and not a configuration of rational curves). 
The remaining cases will be studied in a second part. Notice that $n(-K_X) = 3$ and $K_X^3 = 0$ is the only case when ``abundance'' 
$\kappa (-K_X) = \nu (-K_X)$ does not hold. The surface analogues are $\bP_2$ blown up in 9 points without elliptic fibration
resp.\ elliptic ruled surfaces of very special type.\\
Concerning boundedness (in the rational connected case), we are immediately done by classification if $X$ admits a Mori contraction
$\varphi: X \to Y$ of fiber type, i.e. $\dim Y \leq 2.$ If $\varphi$ is birational, we want to proceed by induction on 
the Picard number. In most cases $-K_Y$ is again nef, so that the induction is no problem. However there are two exceptions, namely
when $X \to Y$ is the blow up of a rational curve $C \subset Y$ with normal bundle $\sO(-1) \oplus \sO(-2)$ resp.\ $\sO(-2) \oplus \sO(-2).$
The first case does not create any difficulty because after some birational transformation it leads to a bounded situation. 
However the $(-2,-2)$-case needs further consideration which will be carried out in the second part of this paper. Therefore
at the moment we obtain boundedness modulo boundedness of threefolds with $n(-K_X) = 3$ resp.\ of those threefolds carrying a
$''(-2,-2)''$-contraction.

\setcounter{tocdepth}{1}
\tableofcontents

\section{Preliminaries}

In \cite{8authors} the following reduction theorem is proved.

\begin{theorem} Let $L$ be a nef line bundle on a normal projective variety $X.$ Then there exists an almost holomorphic dominant
meromorphic map $f: X \rightharpoonup B$ with connected fibers such that 
\begin{enumerate}
\item $L$ is numerically trivial on all compact fibers $F$ of $f$ of dimension $\dim X - \dim B$
\item for a general point $x \in X$ and every irreducible curve $C$ passing through $x$ such that $\dim f(C) > 0,$ we have
$L \cdot C > 0.$ 
\end{enumerate}
The map $f$ is unique up to birational equivalence of $B.$ 
\end{theorem}

Recall that a meromorphic map $f: X \rightharpoonup Y$ is almost holomorphic if there exists an open
non-empty set $U \subset X$ such that $f \vert U$ is holomorphic and proper. In particular $\dim B$
is an invariant of $L$ and we set $n(L) = \dim B,$ the nef dimension of $L$.

\begin{proposition}
Let $p: X \to Y$ be a proper surjective morphism of normal projective varieties and let $L$ be a nef
line bundle on $Y$ with nef reduction $f: Y \rightharpoonup B$. Then the Stein factorization of 
$ f \circ p$ gives a nef reduction for $p^* L$; in particular $n(p^* L) = n(L)$.
\end{proposition}

\proof
Obviously $p^* L$ ist numerically trivial on compact fibers of $ f \circ p$ of the expected dimension.
Let $X \stackrel{g}{\rightharpoonup} A \stackrel{q}{\to} B$ be the Stein factorization of 
(a desingularization of) $f \circ p$. Then $p^* L$ is numerically 
trivial on the general fiber of $g$ so the nef reduction of $p^* L$ must factor via $g$. Let $x \in X$ be
a general point and $C \subset X$ an irreducible curve through $x$ with $\dim g(C) > 0$. As $q$ is finite
$q(g(C))$ is again a curve, so $p(C)$ is a curve which is not contracted by $f$, i.~e.\ $L \cdot p_*(C) >0$.
Now the projection formula implies $p^* L \cdot C >0$ and $g$ is a nef reduction for $p^* L$.
\qed

\begin{corollary}
Let $p: X \to Y$ be an \'etale covering of projective manifolds. Then $n(\pm K_X) = n(\pm K_Y)$.
\end{corollary}

\begin{defprop}
Let $X$ be a projective manifold (or a variety with $\bQ$-factorial canonical singularities say) and 
let $D$ be a nef divisor on $X$. We define the numerical dimension of $D$ to be 
$\nu(D) = \max \, \{ n \, | \, D^n \not\equiv 0 \}$. Then we always have the inequalities 
$\kappa(D) \leq \nu(D) \leq n(D)$ . 
We call $D$ \emph{good} if $\kappa(D) = \nu(D)$, otherwise we call it \emph{bad}.
If $D=\pm K_X$ is good then it is semi-ample, i.~e.\ some multiple is generated by global sections.
\end{defprop}

\proof \cite[2.2 and 6.1]{Ka85}, resp. \cite[2.8]{8authors} for the inequality $\nu (D) \leq n(D).$ 
\qed

By the Abundance Conjecture, the canonical bundle is never expected to be bad, whereas the anticanonical bundle
\emph{can} be bad.
\medskip

The classification of algebraic surfaces with nef anticanonical bundle is of course an easy consequence of the
Kodaira-Enriques classification:

\begin{proposition} Let $X$ be a smooth projective surface with $-K_X$ nef and not numerically trivial. Then the following assertions are equivalent 
\begin{enumerate}
\item $X$ admits an elliptic fibration 
\item $n(-K_X) = 1$ 
\item either after a finite \'etale cover $X \simeq \bP_1 \times A$ with an elliptic curve $A$ or $X$ is $\bP_2$ blown up in 9 points
such that some multiple $-mK_X$ is generated by global sections.
\end{enumerate}
In particular $\kappa(-K_X)=\nu(-K_X)=n(-K_X)=1$ and the nef reduction can be chosen holomorphic, not only almost holomorphic.
\end{proposition}

As a corollary, $n(-K_X) = 2$ if and only if either $-K_X$ is big or $X$ is $\bP_2$ blown up
in 9 points without elliptic fibration or if $X = \bP(E)$ with a semi-stable bundle $E$ over an elliptic curve $A$ which cannot
be written -- after twist -- in the form $\sO \oplus L$ with $L$ torsion or as non-split extension of a trivial line bundle with a line
bundle of degree $1$. Hence:

\begin{proposition} Let $X$ be a smooth projective surface with $-K_X$ bad. Then $\kappa(-K_X)=0$, $\nu(-K_X)=1$, 
$n(-K_X)=2$ and $X$ is one of the following:
\begin{list}{}{}
\item {\bf Case A)} $X$ is $\bP_2$ blown up in 9 points in sufficiently general position (possibly infinitely near) or
\item {\bf Case B)} $X = \bP(E)$, $E$ a rank 2 vector bundle over an elliptic curve which is defined by an extension $0 \to \sO \to E \to L \to 0$ 
with $L$ a line bundle of degree 0 and either
\begin{list}{}{}
\item {\bf B.1)} $L = \sO$ and the extension is non-split or
\item {\bf B.2)} $L$ is not torsion.
\end{list}
\end{list}
\end{proposition}

\begin{remark} In these cases, the structure of the unique element $D$ in $|-K|$ is as follows:
\begin{list}{}{}
\item {\bf Case A)} The image of $D$ in $\bP_2$ is the unique cubic curve containing the 9 points and every point is
a simple point on the cubic. $D$ is either smooth elliptic or a configuration of rational curves.
Every component contains exactly $3d$ points where $d$ is the degree of the component.
\item {\bf B.1)} $D=2C$ and $C$ is smooth elliptic.
\item {\bf B.2)} $D=C_1 + C_2$ where the $C_i$ are smooth elliptic curves which do not meet.
\end{list}
\end{remark}

The remaining case is $-K_X$ big and nef which implies that $X$ is $\bP_2$ blown up in at most 8 points in almost general position, 
$\bP_1 \times \bP_1$ or the Hirzebruch surface $F_2$.

\section{Nef reduction for the anticanonical bundle}

In this section we study the nef reduction of a projective {\it three}fold $X$ with nef anticanonical bundle $-K_X$ and prove that
the reduction map can be taken holomorphic.

\begin{theorem} Let $X$ be a projective threefold with $-K_X$ nef.
Then there exists a
holomorphic map $f: X \to B$ to a normal projective variety $B$
such that 
\begin{enumerate}
\item $-K_X$ is numerically trivial on all fibers of $f$
\item for $x \in X$ general and every irreducible curve $C$ passing through $x$ such that $\dim f(C) > 0,$ we have $-K_X \cdot C > 0.$
\end{enumerate}
In case $X$ is rationally connected and $n(-K_X) = 1$ or $2$ then even some multiple $-mK_X$ is spanned by global sections, so
that we can take $f$ to be (the Stein factorisation of) the map defined by the sections of $-mK_X.$   
\end{theorem}

\begin{definition} Let $X$ be a smooth projective threefold. A $(-2,-2)$-contraction on $X$ is a blow-up $\varphi: X \to Y$ of the smooth 
threefold $Y$ along a smooth rational curve $C$ with normal bundle $N_{C/Y} = \sO(-2) \oplus \sO(-2).$ 
\end{definition} 

\begin{remark} The background of this definition is the following. Let $X$ be a smooth projective threefold, $\phi: X \to Y$ 
the blow-up of a smooth curve $C$ in the projective threefold $Y.$ Then $-K_Y$ is nef unless possibly $\phi$ is a $(-2,-2)$-contraction
or $C = \bP_1$ with normal bundle $\sO(-1) \oplus \sO(-2).$ This last case is usually easy to be dealt with. 
\end{remark}

\proof If $X$ is not rationally connected, then the assertions follow from direct classification,
see sect.\ 3 and 4. So we will assume $X$ to be rationally connected -- {\it actually we shall use rational connectedness only if $K_X^2 = 0.$} 
We start with an almost holomorphic nef reduction $f: X \rightharpoonup B$ over the normal projective variety $B$ of
dimension $b.$ If $b = 1,$ then $f$ is automatically holomorphic,the spannedness been proved in (5.2), so we may assume $b = 2.$
Consider the general fiber $C$ of $f$, an elliptic curve, and form the associated family $(C_t)_{t \in T}.$ To be precise, we consider the
graph of this family $$q: \sC \to T$$ with projection
$ p: \sC \to X.$ Then $p(q^{-1}(t)) $ is a compact fiber $C_t$ of $f$ for general $t$ and of course 
$p : q^{-1}(t)  \to C_t$ is an isomorphism for all $t.$ After a base change we may assume $T$ smooth
and also we may assume $\sC$ normal. \\
Since $p^*(-K_X)$ is $q$-nef, we easily find a line bundle $\tilde L$ on $T$ and a positive integer $m$ such
that $$ p^*(-mK_X) = q^*(\tilde L).$$
Indeed, we let $L = q_*(p^*(-K_X))^{**};$ notice here that $p^*(-K_X)$ is trivial on the general fiber
of $q$ since $K_X \vert C_t $ is trivial and since $p$ is an isomorphism near the general $C_t.$  
Thus 
$q_*(p^*(-K_X))^{**}$ is really a line bundle on $T$ and we obtain an injection $q^*(L) \to p^*(-K_X).$ This yields a decomposition
$$ p^*(-K_X) = q^*(L) + D $$
with $D$ effective coming from multiple fibers. Hence $mD = q^*(D')$ and we put $\tilde L = mL + D'.$ 
\\
If now $K_X^2 \ne 0,$ then $L^2 > 0$ and therefore $\kappa (L) = 2. $ This implies $\kappa (-K_X) = 2$ (since $p$ has degree 1)
so that the numerical Iitaka dimension and the Iitaka dimension coincide. Thus $-K_X$ is {\it good}.
By [Ka85,6.1], $-K_X$ is therefore semi-ample, i.e. some $-mK_X$ is spanned. 
\smallskip

It remains to treat the case $K_X^2 = 0.$
Here $L^2 = 0,$ but of course we cannot say that
$\kappa (L) = 1.$ If however we know that $\kappa (-K_X) = 1,$ then the same arguments as above show that
$-K_X$ is semi-ample. To get more informations, we consider a Mori contraction $\varphi: X \to Y.$ 
Since $K_X^2 = 0, $ we rule out $\dim Y = 1$ and also $\varphi$ cannot contract a divisor to a point.
In other words, $\varphi$ is a conic bundle or the blow-up of a smooth curve $C$ in the smooth 3-fold
$Y.$ \vskip .2cm \noindent
(1) Suppose first that $\varphi $ is  birational. Then $-K_Y$ is again nef unless $\varphi$ is a $(-2,-2)$-contraction since 
by $K_X^2 = 0,$ the second exception in (2.2) cannot occur. We also note that $K_Y^2 = [C].$ 
If $-K_Y$ is nef, the Kawamata-Viehweg vanishing theorem gives
$$H^2(-K_Y) = 0.$$
On the other hand, $\varphi_*(-K_X) = \sI_C \otimes -K_Y,$ hence by virtue of $R^q\varphi_*(-K_X) = 0$,
we obtain the exact sequence
$$ H^1(C,-K_Y \vert C) \to H^2(-K_X) \to H^2(-K_Y).$$
Since $K_Y \cdot C = 2 - 2g(C) = \deg (-K_C) $ (see the proof of 5.3 for the detailed computations),
we obtain $ h^1(-K_Y \vert C) = h^0(K_Y \vert C + K_C) \leq 1.$ 
Now Riemann-Roch gives
$$ \chi(-K_X) = 3.$$
Here we used the rational connectedness to obtain $\chi(\sO_X) = 1.$ 
Putting things together, we obtain $h^0(-K_X) \geq 2$ and therefore $\kappa (-K_X) = 1.$ In case
of a $(-2,-2)$-contraction, we verify the vanishing $H^2(-K_Y) = 0$ by hand, then all arguments are the same. 
By duality, we check $H^1(2K_Y) = 0.$ Let $H$ be a general hyperplane section. Then $-2K_Y+H$ is ample, at least after
substituting $H$ by a multiple; hence $H^1(2K_Y-H) = 0.$ Since $H$ does not contain the curve $C$ and since $K_Y^2 = C,$ 
the restriction $-K_Y \vert H$ is big and nef, hence $H^1(2K_Y \vert H ) = 0.$ Thus $H^1(2K_Y) = 0.$ 
\vskip .2cm \noindent
(2) Now consider the case of a conic bundle $\varphi: X \to Y$ with discriminant locus $\Delta \subset Y$. As $X$ is
rationally connected, $Y$ must be rational. By the formula $\varphi_* (K_X^2) \equiv -(4K_Y + \Delta)$ we know that
$|-4K_Y|$ contains the reduced element $\Delta$, hence $-K_Y$ must be nef (cf.\ Lemma 2.5). So $Y$
is either $\bP_1 \times \bP_1$, $F_2$ or $\bP_2$ blown up in at most 9 points in almost general position.
As above, Riemann Roch gives $\chi(-K_X)=3$. If $h^0(-K_X) \geq 2$, then $\kappa(-K_X) = \nu(-K_X) = 1$ and
we are done. So we have to rule out the two possibilities
\smallskip

{\bf A)} $h^2(-K_X) \geq 3$ resp.
\smallskip

{\bf B)} $h^0(-K_X)=1$ and $h^2(-K_X)=2$.
\smallskip

\noindent We consider the rank 3 vector bundle $V = \varphi_*(-K_X)(-K_Y)$. Using duality and the vanishing of the
higher direct image sheaves $R^i\varphi_*(-K_X)$ we calculate $h^0(V^*) = h^2(\varphi_*(-K_X)) = h^2(-K_X)$ so we know in case
{\bf A)} that $V^*$ has at least 3 sections and in case {\bf B)} that $V^*$ has 2 sections.\\
As $-K_X$ restricted to a fiber of $\varphi$ is $\sO(1)$ on the conic we can recover $X$ as a hypersurface in $\bP(V)$. By
construction $X$ is linear equivalent to $2 \xi + \pi^*D$ for some divisor $D$ on $Y$ and
$(\xi + \pi^*K_Y)_{|X} = -K_X$ which determines $D$: Using the canonical bundle formula we calculate
$$ -K_X = (- K_{\bP(V)} - X)_{|X} = (3 \xi - \pi^*(K_Y + \det V) - 2\xi -\pi^*D)_{|X}$$
$$ = (\xi - \pi^*(K_Y + \det V + D)_{|X}$$
hence $ D = -2 K_Y - \det V$. Second, the condition $K_X^2=0$ reads as
$$ 0 = (\xi + \pi^*K_Y)^2 \cdot (2 \xi + \pi^*D)$$
$$ = \xi^2 \cdot \pi^*(2 \det V + 4K_Y + D) + \xi \cdot \pi^*(-2 c_2 + 2K_Y^2 + 2K_Y \cdot D)$$
which implies $D = -2 \det V - 4K_Y = 2 D$ i.~e.\ $D=0$ as well as $c_2 = K_Y^2$.\\
In case {\bf A)} $V^*$ has 3 sections $s_1, s_2, s_3$. We consider a general line $C \subset Y$. Lemma 2.6 implies that
$V$ restricted to $C$ is nef hence of type $\sO(a) \oplus \sO(b) \sO(c)$ with nonnegative integers $a,b,c$.
As $V^*$ has 3 sections we conclude $a=b=c=0$ hence $\det V \cdot C =0$. But then $K_Y \cdot C = 0$ using our
numerical condition which is impossible.\\
In case {\bf B)} we still know that $V^*$ has two sections $s_1, s_2$. Again $V$ restricted to a general $\bP_1$
is nef and therefore $V$ splits as $\sO \oplus \sO \oplus \sO(a)$ with $a \geq 0$. In particular
$s_1 \wedge s_2$ does not vanish identically. This means that the cokernel $\sL$ in the sequence
$$ 0 \to \sO \oplus \sO \stackrel{(s_1,s_2)}{\longrightarrow} V^* \to \sL \to 0$$
has generic rank 1. From the sequence and our numerical conditions above we calculate $c_1 (\sL) = -c_1 (V) = 2K_Y$ 
and $c_2(\sL) = c_2(V) = K_Y^2$. Dualising we obtain an injection $ 0 \to \sL^* \to V$
with $\sL^*$ locally free and in fact $\sL^* = -2K_Y$ because on a rational surface numerical and linear equivalence coincide.
Twisting by $K_Y$ we finally obtain an injection
$$ 0 \to \sO(-K_Y) \to \varphi_*(-K_X).$$
If $-K_X$ has more that one section $\kappa(-K_X)=1$ and we are done. Therefore $h^0(-K_Y)=1$. An application of Riemann-Roch
shows that this is only possible if $K_Y^2=0$ hence $c_2(V)=0$ and already $\sL$ is locally free. So we get:
$$ 0 \to \sO(-K_Y) \to \varphi_*(-K_X) \to \sO(K_Y) \oplus \sO(K_Y) \to 0$$
Now we have two possibilities: Either this splits i.~e.\ $\varphi_*(-K_X) = \sO(-K_Y) \oplus \sO(K_Y)^{\oplus 2}$ and
$-2K_X$ has at least 3 sections and again we are done. Or the sequence doesn't split which is only possible if
$h^1(-2K_Y)=1$ which is equivalent to $h^0(-2K_Y)=2$ using Riemann-Roch. Writing down the usual filtration for $S^2(\varphi_*(-K_X))$
we obtain
$$ 0 \to F^1 \to S^2(\varphi_*(-K_X)) \to \sO(2K_Y)^{\oplus 3} \to 0$$
$$ 0 \to \sO(-2K_Y) \to F^1 \to \sO^{\oplus 2} \to 0$$
and $\kappa(-K_X)=1$ once again.
\qed

In the proof of Theorem 2.1 we have actually shown the following

\begin{corollary} Let $X$ be a smooth projective threefold with $-K_X$ nef. Suppose $K_X^2 = 0$
and $X$ rationally connected. Then $\nu(-K_X) = \kappa (-K_X) = 1$ and therefore $-K_X$ is semi-ample,
inducing an abelian or K3-fibration over $\bP_1.$ 
\end{corollary} 

During the proof of Theorem 2.1 we used the fact that for a conic bundle $X \to Y$ with $K_X^2=0$ we know
that $-K_Y$ is nef. This follows from the following lemma:

\begin{lemma}
Let $Y$ be a smooth projective surface and assume that there exists a reduced divisor $\Delta \equiv -mK_Y$ for 
some $m \geq 4$. Then $-K_Y$ is nef.
\end{lemma}

\proof Let $C \subset Y$ be a curve which is contained in $\Delta$, i.~e.\ $\Delta = C + \Delta'$ for some
effective divisor $\Delta'$ which does not contain $C$. Then $(1-m) K_Y \cdot C = K_Y \cdot C + \Delta \cdot C = 
\deg K_C + \Delta' \cdot C \geq -2$. Therefore $-K_Y \cdot C \geq -2 / (m-1)$ and the integer $-K_Y \cdot C$
is nonnegative.
\qed

Another ingredient of the proof is the following specialization of $\cite[3.21]{DPS94}$.

\begin{lemma} Let $\varphi: X \to Y$ be a Mori contraction which is a conic bundle, $X$ a smooth projective threefold 
with $-K_X$ nef and let $V$ be the vector bundle $\varphi_*(-K_X)(-K_Y)$. Furthermore let $\bP_1 \cong C \subset Y$ such that $X_C = \varphi^{-1}(C)$
is smooth. Then $V_{|C}$ is generated by global sections, in particular it is nef.

\end{lemma}

\proof We consider the induced conic bundle $\varphi_C: X_C \to C$. Let $l = \varphi_C^{-1}(y)$ be any fiber. 
We want to show that $(\varphi_C)_*(-K_{X|C})(-K_{Y|C})$ is generated by its global sections, i.~e.\ 
that every section of $( -K_X \otimes \varphi^* (-K_Y) )_{|l}$ lifts to $X_C$. We will show the vanishing of
$H^1(X_C, -K_{X|X_C} - \varphi_C^*(-K_{Y|C}) -l)$ which gives the desired extension property. If we write
$$-K_{X|X_C} - \varphi_C^*(-K_{Y|C}) -l) = K_{X_C} + L$$
then some easy calculation gives 
$$ L = -2 K_{X|X_C} + \varphi_C^*(-K_C -y)$$
As $-K_{X}$ is nef and $\varphi$-ample $L$ is ample and Kodaira vanishing gives $H^1(K_{X_C} + L)=0$.
\qed

\begin{remark}
We used several times the fact that the restriction of $\xi = \sO_{\bP(V)}(1)$ to $X$ is $-K_X$. 
If $l$ is a conic then $-K_X \cdot l = \xi \cdot l = 2$ so $\xi_{|X}$ and $-K_X$ only differ by $\varphi^* M$ for some $M$.
Consider the relative Euler sequence
$$ 0 \to \Omega^1_{\bP(V)/Y}(1) \to \pi^* V \to \sO_{\bP(V)} (1) \to 0$$
Since $H^i( l, \Omega^1_{\bP_2}(1)_{|l}) =0$ for $i=0,1$ which follows from the Bott formula we know that
$\varphi_*( \Omega^1_{\bP(V)/Y}(1)_{|X}) = R^1\varphi_*( \Omega^1_{\bP(V)/Y}(1)_{|X}) = 0$.
So if we restrict the sequence to $X$ and push it down to $Y$ we obtain an isomorphism
$$\varphi_*(-K_X) = V = \varphi_* \varphi^* V = \varphi_* (\pi^*V_{|X}) \stackrel{\cong}{\to} \varphi_*(\xi_{|X})$$
which implies that $M$ is torsion hence trivial as $Y$ is simply connected.
\end{remark}

\section{Non-rationally connected threefolds with vanishing irregularity} 

\begin{theorem} Let $X$ be a smooth projective threefold with $\tilde q(X) = 0.$ Suppose that $-K_X$ is nef but $K_X \not \equiv 0.$ 
Then either $X$ is rationally connected or after finite \'etale cover, $X \simeq \bP_1  \times S$ with $S$ a K3 surface. 
\end{theorem}

\proof  Since $K_X \not \equiv 0,$ its Kodaira dimension $\kappa (X) = - \infty,$ hence $X$ is uniruled. Let $(C_t)$ be a covering family 
of rational curves, providing an almost holomorphic quotient $f: X \rightharpoonup S$ to a smooth variety $S$ of dimension 1 or 2. 
If $\dim S = 1,$ then then $S \simeq \bP_1$ by $\tilde q(X) = 0.$ Since the fibers of $f$ are rationally connected, we conclude that 
$X$ is rationally connected. So let $\dim S = 2;$ we may assume $S$  minimal. Since $\tilde q(X) = 0, $ also $\tilde q(S) = 0$.
Moreover by Zhang [Zh96], $\kappa (S) \leq 0,$ see [DPS01,4.12]. If $\kappa (S) = - \infty,$ then $S$ must be rational, and hence $X$ is rationally
connected. Thus we are reduced to $S$ being K3 or an Enriques surface. After a finite \'etale cover we may assume $S$ to be K3. \\
By Mori theory (in the relative version for $f$) we can find a sequence $g: X \rightharpoonup X'$ of birational contractions and
flips such that $X'$ has a $\bP_1$-fibration $h: X' \to S,$ which is the contraction of an extremal ray. 
By [PS98,2.1,2.2] $-K_{X'}$ is almost nef; i.e. $-K_{X'} $ is nef except for finitely many rational curves. Hence the arguments of
[PS98,1.9,1.10] apply and show that $X' \to S$ is actually a $\bP_1$-bundle. In particular, $X'$ is smooth. Then we can write
$$ X' = \bP(E) $$
with a rank 2 vector bundle $E$ over $S.$ Now the almost nefness of $-K_X$ implies that $S^2E \otimes \det E^* $ is nef on all curves
except for finitely many rational curves $C_i  \subset S.$ Using $\bQ$-bundles we can say that
$$ E_0 := E \otimes {{\det E^*} \over {2}} $$
is almost nef. Since moreover $-K_X$ is pseudo-effective, also $-K_{X'} $ is pseudo-effective, so that $\sO_{\bP(E_0)}(1) $ is
pseudo-effective. Then by [DPS01,6.7(c)], $E_0$ is numerically flat.  \\
NB. It is clear that we may argue on the level of $\bQ$-bundles; alternatively note that, fixing $A$ ample on $S,$ then
$$ H^0(S^m(S^2E \otimes \det E^*) \otimes A)) \ne 0 $$
for large $m$, hence $\sO_{\bP(S^2E \otimes \det E^*)}(1) $ is pseudo-effective and also almost nef in the sense of [DPS01], so that
[DPS01,6.7(c)] applies to give the numerical flatness of $S^2E \otimes \det E^*.$ \\
Now, $S$ being simply connected, it follows that $S^2E \otimes \det E^* = \sO_S^3,$ and we claim that then $\bP(E) \simeq \bP_1 \times S$ : 
consider a smooth member $D \in \vert -K_{X'} \vert, $ given by a general section in $S^2E \otimes \det E^*.$ Now $D$ maps $2:1$ 
onto $S$ and $K_D = \sO_D$, so $D \to S$ is \'etale. Since $\pi_1(S) = 0,$ $D$ must be disconnected with $2$ components, hence
$-K_{X'}$ is divisible by $2$ and therefore $E \otimes \det E^*$ exists as a vector bundle and is trivial. \\
Alternatively, the three sections of $-K_{X'}$ give a map $X' \to \bP_2$. Since $-K_{X'}^2 = 0,$ the image of this map must be $ \bP_1$
and we conclude. \\
Finally we show that $X = X'.$ So let $g_m: X_m \rightharpoonup X'$ be the last contraction of $g.$ We know again that $-K_{X_m} $ is almost nef.
This leads immediately to a contradiction by considering a surface $S_x = \{x\} \times S$ such that $S_x$ meets the center of $g_m$ 
and by computing canonical bundles. 
\qed

\begin{corollary}  Let $X$ be a smooth projective threefold with $\tilde q(X) = 0.$ Suppose that $-K_X$ is nef but $K_X \not \equiv 0.$ 
If $X$ is not rationally connected, then either $X = \bP_1 \times S$ with
$S$ a K3 or an Enriques surface or $X$ is a non-trivial $\bP_1$-bundle over an Enriques surface $S$ which
is trivialized by the universal 2:1-cover $\tilde S \to S$. In all cases we have $n(-K_X) = n(-K_{\tilde X}) = 1$. 
\end{corollary}

\proof Let $\tilde X \to X$ be the universal cover; then $\tilde X = \bP_1 \times \tilde S$ with a
K3 surface $\tilde S$ by (3.1). Let $\varphi: X \to S$ be a Mori contraction; then $\varphi $ lifts to
$\tilde X$ and therefore must be a $\bP_1$-bundle. Moreover we obtain an \'etale cover $\tilde S \to S.$ 
So $S$ is K3 or Enriques and we are in the situation as described in the corollary.
\qed

\section{Threefolds with positive irregularity} 

\begin{setup} {\rm
Here we consider smooth projective threefolds $X$ with $-K_X$ nef such that $\tilde q(X) > 0.$ After passing to a finite \'etale cover
we may and will assume that $q(X) = \tilde q(X).$ We let $\alpha: X \to A$ denote the Albanese, which is a surjective submersion by
[PS98]. In particular $q \leq 3$ and $q = 3$ iff $X$ is abelian. So we need only to consider the cases $q = 1$ and $q = 2.$ 
Additionally we consider the invariant $n(-K_X).$ If $n(-K_X) = 0, $ then $K_X \equiv 0,$ so we suppose $n(-K_X) > 0.$ Then we have
a non-trivial nef reduction $f: X \rightharpoonup S$, at least if $n(-K_X) \ne 3$ (in which case $S = X).$ }
\end{setup} 
 
\begin{theorem} Let $X$ be a smooth projective threefold with $-K_X$ nef. Suppose
$$1 \leq q = \tilde q \leq 2$$ and
$$0 < n(-K_X) < 3.$$
Then the nef reduction $f: X \to S$ can be taken holomorphic and after a finite \'etale cover, $X$ is one of the following
\begin{enumerate}
\item $q(X) = 1,$ $n(-K_X) = 2,$ and $X \simeq A \times S$ ($A$ elliptic, $S$ rational with $n(-K_S)= 2$, i.e. $\bP_2,$ $\bP_1 \times \bP_1$,
del Pezzo or $-K_S$ nef without elliptic fibration).
\item $q(X) = 1,$ $n(-K_X) = 1,$ and $X \simeq A \times F$ where $F$ is $\bP_2$ blown up in 9 points such that $-K_F$ is nef and $F$
admits an elliptic fibration. Here $S$ is the image of the elliptic fibration on $F.$
\item $q(X) = 2,$ $n(-K_X) = 2.$ Then $X = \bP(E) \times B_2$ with $B_i$ elliptic curves and $E$ a numerically flat bundle over $B_1$
with $n(-K_{\bP(E)}) = 2$ (i.e. $E$ is non-trivial even after finite \'etale cover). Here $S = \bP(E).$ 
\item $q(X) = 2,$ $n(-K_X) = 1.$ Then $X = A \times \bP_1 $ ($A$ abelian), and $S = \bP_1.$
\end{enumerate}
\end{theorem}

\proof We consider the almost holomorphic nef reduction $f: X \rightharpoonup S$ to a smooth curve or 
a normal surface. \\
\\
{\bf Case I:} $q = 1.$ \\ 
Then all fibers $F$ of $\alpha$ are smooth rational surfaces with $-K_F$ nef and in particular $K_F^2 \geq 0.$ \\
{\it Subcase I.1:} $F = \bP_2.$ So $\alpha$ is a $\bP_2$-bundle; we write $X = \bP(E)$ with a rank 3-bundle $E$ over $A.$ 
Then $E \otimes {{\det E^*} \over {3}}$ is nef with $c_1 = 0,$ hence numerically flat. 
If $n(-K_X) = 2,$ then $\dim S = 2$ and consider a general fiber $C$ of $f.$ Since $f$ is holomorphic with $K_X \cdot C = 0,$
$C$ is a smooth elliptic curve and therefore an \'etale multi-section of $\alpha,$ hence a section after a finite \'etale
cover of $A.$ Then we obtain a 2-dimensional family of disjoint sections, and hence $X = A \times S$ \\
The case $n(-K_X) = 1$ is obviously impossible by dimension reasons since $-K_X$ is $\alpha$-ample and trivial on the general fiber of $f.$ \\
{\it Subcase I.2:} $F = \bP_1 \times \bP_1.$  Then we have a factorisation (see [PS98])
$$ X {\buildrel {g} \over {\la }} W {\buildrel {h} \over {\la}} A $$
with $g$ and $h$ both $\bP_1$-bundles. Then $X = \bP(E) $ with a rank 2-bundle $E$ over $W$; moreover $-K_W$ is nef. Since $-K_X$ is $\alpha$-ample,
it is clear that (as in case I.1) $\dim f(F_{\alpha}) =  2,$ hence $n(-K_X) = 2.$ Let again $C$ be a general fiber of $f,$ a smooth
elliptic curve. Then after finite \'etale cover of $A$, $C$ is a section of $\alpha,$ hence already $W = \bP_1 \times A.$ 
These elliptic curves define an irreducible family $(C_t)_{t \in T}$ with $T$ compact, a priori only the general $C_t$ is a smooth
elliptic curve. However $-K_X \cdot C_t = 0$ for all $t;$ hence $C_t$ cannot contain a component in a fiber $F$ of $\alpha.$ On the other
hand, $C_t  \cdot F = 1,$ hence every $C_t$ is a section of $\alpha.$ Now consider the family $(g(C_t)).$ This is a complete family of
elliptic curves, i.e. with no degeneracies. Therefore $W = \bP_1 \times A$ and $g(C_t)$ is a fiber of the projection to $\bP_1.$ 
Also $\bP(E \vert c \times A) = \bP_1 \times A,$ hence after renormalizing, $E \vert c  \times A$ is trivial for all $c$ and so 
$E = p_1^*(E') $ with a vector bundle $E'$ over $\bP_1.$ Hence $X = \bP(E') \times A$, and consequently $\bP(E') = \bP_1 \times
\bP_1.$
Therefore $X = \bP_1 \times \bP_1 \times A.$ \\
{\it Subcase I.3:} $F$ is del Pezzo with $K_F^2 \leq 7.$ Then we have a factorization  
$$ X {\buildrel {g} \over {\la }} W {\buildrel {h} \over {\la}} A $$
with $g$ the blow-up of some multi-sections and $h$ a $\bP_2$-bundle. By the same reason as before, $n(-K_X) = 2.$ Again we get a lot of
multi-sections of $h$, which gives $W = \bP_2 \times A$ (after \'etale cover), hence $X = \bP_2(x_1, \ldots, x_r) \times A.$ \\
{\it Subcase I.4} The same arguments still work if $-K_F$ is just nef as long as $\dim S = 2,$ i.e. $n(-K_X) = 2.$ 
So suppose $\dim S = 1.$ Note that this is only possible if $-K_F$ is not big, i.e. $K_F^2 = 0$ and if $F$ admits an elliptic
fibration to $S = \bP_1.$ 
We still have a factorization  
$$ X {\buildrel {g} \over {\la }} W {\buildrel {h} \over {\la}} A $$
with $g$ the blow-up of some multi-sections and $h$ a $\bP_2$-bundle. 
Consider the nef reduction $f: X \to S \simeq \bP_1.$ Then the general fiber $F$ is a smooth surface with $K_F \equiv 0.$
Since $F$ projects onto the elliptic curve $A,$ the surface $F$ must by hyperelliptic or a torus and actually it is a product after
finite \'etale cover. Let $E$ be the exceptional locus of $g.$ Then $E \cap F$ is a union of multi-sections of $F \to A$ and we are going to determine
its structure. \\
Consider the last blow-up $g_r : X = X_r \to X_{r-1},$ blowing up the \'etale multi-section $C_r$ of $X_{r-1} \to A.$ Let $E_r$ be the corresponding divisor
in $X.$ Since $f(E_r) = S,$ we have an \'etale cover $E_r \to S \times A$ given by $(f \vert E_r) \times (h \circ g \vert E_r)$. 
Hence $E_r \cap F$ is an elliptic curve, an \'etale multi-section of $F \to A.$ Therefore we find a 1-dimensional (non-complete) family
of disjoint multi-sections of $F \to A$ not meeting $E \cap F.$ Varying $F$ and proceeding by induction on the blow-ups belonging to $g,$ 
we obtain a 2-dimensional family of disjoint \'etale multi-sections of $W \to A$ (not meeting the exceptional locus of $g$ in $W$).
Hence $W$ is a product after finite \'etale cover and we want to conclude that then $X$ is already a product (up to finite
étale cover). This is clear if we always blow up curves of type $A_p = A \times \{p\}$. So assume this is not the case 
i.~e.\ we blow up a curve $C_i$ which is not of this type. Then we may find some curve $A_p$ such that $C_i \cap A_p$ is finite.
But if $\hat{A}_p$ denotes the strict transform we then calculate
$$ -K_{X_i} \cdot \hat{A}_p = ( g_i^* (-K_{X_{i-1}}) - E_i) \cdot \hat{A}_p = - K_{X_{i-1}} \cdot A_p - E_i \cdot \hat{A}_p = - E_i \cdot \hat{A}_p <0$$
(here $-K_{X_{i-1}} \cdot A_p = 0$ because inductively we may assume that $X_{i-1}$ is a product). Now this contradicts $-K_{X_i}$ nef.\\
{\bf Case II:} $q = 2$. Now $\alpha: X \la A$ is a $\bP_1$-bundle. After a finite \'etale cover, we can write $X = \bP(E)$
with a numerically flat rank 2-bundle $E$. \\
{\it Subcase I:} \  $n(-K_X) = 1.$ In that case we have again a holomorphic map $f: X \la S = \bP_1$. The general fiber $F_f$ is an 
\'etale cover of $A,$ so after another \'etale base change, general fibers of $f$ and $\alpha$ meet transversally at one point.
In other words, $\alpha$ has many disjoint sections and thus $E$ is trivial (after a twist). So $X = S \times A =  \bP_1 \times A.$\\
{\it Subcase II:} \ $n(-K_X) = 2.$ Then $\alpha$ has a 2-dimensional family of sections so that $A$ carries a famliy of elliptic curves.
Now by Poincar\'e reducibility, $A = B_1 \times B_2$ is a product of elliptic curves $B_i;$ possibly after finite \'etale cover.
Then we argue similarly as in Subcase I.2 to obtain the product structure $X = \bP(F) \times B_2$ with $F$ a semi-stable bundle on
$B_1 $ (such that $E = p_1^*(F)).$ 
\qed

\begin{re} {\rm In the case $n(-K_X) = 3$ we cannot expect such precise results. Here $-K_X$ is positive on all covering families of
generically irreducible curves.  
Let us consider e.g. the situation that
$q(X) = 2.$ With the notations as before and after a finite \'etale cover, $X = \bP(E)$ with a rank 2-bundle over the Albanese torus
$A$ such that $E$ is nef with $c_1(E) = 0 $ (see [CP91]). So $E$ is numerically flat. Fix a curve $C \subset A.$ Then there exists a
moving curve $B \subset \bP(E \vert C)$ with $K_X \cdot B = 0$ if and only if after normalizing $C$ and after a finite \'etale cover
of the normalization, the bundle $E \vert C$ splits. Let us say that $E$ is almost trivial on $C.$ Then we obtain: \\
$n(-K_X) = 3$ if and only if there is at most a countable number of curves $C \subset A$ such that $E \vert C$ is almost trivial. }
\end{re}

\begin{corollary} Let $X$ be a smooth projective threefold with $-K_X$ nef and $K_X \not \equiv 0.$ 
If $\tilde q > 0,$ then $q > 0.$ The nef reduction is holomorphic, and if $q < \tilde q$ then  
$X$ is a $\bP_1$-bundle over a hyperelliptic surface.  
\end{corollary} 
 
\proof Since $K_X \not \equiv 0,$ we have $\tilde q \leq 2.$ Let $\tilde X \to X$ be a finite \'etale 
cover such that $q(\tilde X) = \tilde q.$ 
From the classification in (4.2) we deduce $\chi(\sO_{\tilde X}) = 0.$ If $q(X) = 0$ then $\chi(\sO_X) 
\geq 1,$ contradicting $\chi(\sO_X) = \frac{1}{m} \chi(\sO_{\tilde X}).$ So we only have to investigate threefolds 
$X$ with $q = 1$ and $\tilde q = 2.$ Looking at the classification, we see that then the Albanese 
map $\alpha: X \to A$ factors as $g \circ f$ with a $\bP_1$-bundle $f: X \to A'$ and an elliptic 
bundle $g: A' \to A$ where $q(A') = 1.$ So $A'$ is hyperelliptic (since $\kappa (A') \leq 0$). The fibers $F$ of
$\alpha$ are $\bP_1$-bundles over elliptic curves with $-K_F$ nef. If $\dim S = 2,$ $X$ carries a 2-dimensional family of elliptic curves.
Then $F$ is necessarily a product, and $S = \bP(V) $ with a rank 2-bundle over $A$ while $X = \bP(g^*(V)).$ The bundle $V$ is moreover non-trivial,
even after twists. If $\dim S = 1,$ then $S = \bP_1$ and $X = S \times A'.$ In all other cases, $\dim S = 0.$ 
\qed

\begin{theorem} Smooth projective threefolds $X$ with $-K_X $ nef, $K_X \ne 0$ and $\tilde q > 0$
form a bounded family up to finite \'etale cover.
\end{theorem}

\proof By virtue of (4.1) we are a priori reduced to the case $n(-K_X) = 3,$ however we will not use this
information. Let $\alpha: X \to A $ be the smooth surjective Albanese map. \\
First assume that $q = 2.$
As already noticed, after possibly a finite \'etale cover, $X$ is of the form $\bP(E)$ with a 
numerically flat rank-2 bundle $E$ over $A$. In particular $c_1(E) = 0$ and $E$ is semi-stable.
This gives boundedness. \\
So from now on we assume $q = 1.$ If $\alpha$ is a $\bP_2$-bundle, the same argument applies. In all
other cases we have the following picture [PS98]. There exists a $\bP_1$-bundle $p: W \to A$ with
$-K_W$ nef and a rank 2-bundle $E$ over $W$ such that the 3-fold $X' = \bP(E)$ has nef anticanonical
bundle and such that $X$ arises from $X'$ by blowing up some \'etale multisection of $\alpha$ (including
the case $X = X'$).
The surfaces with $-K$ nef are bounded, so we may fix $W.$ Next we have to bound $E$ up to twists by
line bundles. Let $F$ be a fiber of $p$ and $C_0$ a section with $C_0^2 $ minimal. We normalize $E$
such that $$ 0 \leq c_1(E \vert F) \leq 1 $$
and $$ 0 \leq c_1(E \vert C_0).$$
Let $e = -C_0^2. $ Since $-K_W$ is nef, we have $-1 \leq e \leq 0.$ If $e = -1,$ then there exists
an \'etale cover of degree $2$ which has $e = 0.$ So we can restrict to $e = 0.$ 
Writing $$ c_1(E) = aC_0 + bF,$$
we have $a = c_1(E \vert F) $ and $b = c_1(E \vert C_0),$ so $0 \leq a,b \leq 1.$ 
In particular $c_1(E)^2 = 2ab = 0$ or $2.$ The fact that $-K_{X'} $ is nef is translated into the nefness
of
$$ E \otimes {{\det E^*} \over {2}} \otimes {{-K_W} \over {2}}.$$
Using $K_W^2 = 0$ and the inequalities $c_1^2 \geq c_2$ and $c_2 \geq 0$ for a nef bundle, the equality
$c_2(E \otimes {{\det E^*} \over {2}}) = 0$ is established. This means
$$ c_1^2(E) = 4c_2(E), $$
hence $c_1^2 (E) = c_2(E) = 0$ and $ab = 0.$ 
If now $E$ is semi-stable with respect to the ample divisor $H = C_0+F$, then we obtain boundedness.
So suppose $E$ is unstable with respect to $H.$ Let $S$ be the maximal destablising subsheaf, a line
bundle. Then we have an exact sequence
$$ 0 \to S \to E \to \sI_Z \otimes Q \to 0 \eqno (S)$$
with a finite set $Z$ and a line bundle $Q$. 
The destabilising property gives $S \cdot H \geq 1.$ Since $K_W \cdot C_0 = 0$, 
$$E \otimes {{\det E^*} \over {2}} \vert C_0$$ 
is numerically flat, hence $E \vert C_0$ is semi-stable. Thus 
$$S \cdot C_0 \leq 
{{1} \over {2}} c_1(E) \cdot C_0 = {{b} \over {2}}, $$ 
so $S \cdot C_0 = 0$ and therefore $S \cdot F 
\geq 1.$ \\
The nefness of $$E \otimes {{\det E^*} \over {2}} \otimes {{-K_W} \over {2}} \vert F$$ yields
$ E \vert F = \sO \oplus \sO $ or $\sO(1) \oplus \sO(-1)$ in case $a = 0$ and
$ E \vert F = \sO \oplus \sO(1) $ in case $a = 1.$ Hence $S \cdot F = 1$  and consequently we have
$S \equiv C_0.$ After tensoring with a topological trivial line bundle we have $S = C_0.$ 
If $a = b = 0,$ then $Q \equiv -C_0.$ Since $c_2(E) = 0,$ $Z = \emptyset $ and now the sequence
(S) proves the boundedness. 
The cases $a = 1, b = 0$ and $a = 0, b = 1$ are done in the same way. \\
Finally we have to deal with the  multi-sections (of degree at most $9$) to be successively
blown up. After some \'etale cover of $A$, the first multisection $C$ is a section of $X' \to A$ and we have to bound $C.$ Let $f:  Z \to X'$ be the
blow up of $C$ and $D = f^{-1}(C).$ 
The equation  
$$ 0 = K_Z^3 = (K_{X'} + D)^3 $$
together with $K_{X'}^3 = 0$ yields $c_1(N_C) = - 3 (K_{X'} \cdot C).$ On the other hand, $c_1(N_C) = - K_{X'} \cdot C $ since $C$
is an elliptic curve. Thus $K_{X'} \cdot C = c_1(N_C) = 0$.  
The nefness of $-K_Z \vert D = -K_D + N_D$ leads easily 
to the statement that $N_{C/X'} $ is numerically flat. 
Let $\tilde C = p(C) \subset W.$ 
Then $\tilde C$ is a section of $W \to A.$ Let $\tilde e$ be the invariant of $E \vert \tilde C.$ 
Let $c = \tilde C^2.$ Then the nefness of $-K_{X'} \vert \bP(E \vert \tilde C)$ is translated into 
$c \geq \tilde e$ if $e \geq 0$ and into $c \geq 0$ if $\tilde e < 0,$ i.e. $\tilde e = -1.$
Let $C_0 $ be a section of minimal self-intersection in $\bP(E \vert \tilde C)$.
Since $-K_{X'}$ is nef and since $-K_{X'} \cdot C = 0,$ we must have $C \equiv C_0.$ 
Since $-K_{X'} \vert \bP(E \vert \tilde C) \equiv 2C_0 + (e+c)F,$ we conclude that $e = -c.$ 
Therefore in total $c \geq 0$ if $\tilde e \geq 0,$ hence $c = \tilde e = 0.$ Hence $\tilde e =
0,-1.$ 
Hence $\tilde C^2 = 0,1$ which proves boundedness of $\tilde C$ and hence of $C.$ \\
The other centers to be blown up are treated in the same way; we leave the details to the
reader.

\qed

\section{Rationally connected threefolds I}

In this section we investigate rationally connected threefolds $X$ with $n(-K_X) = 1;$ they can be viewed as the 3-dimensional 
analogues of the surfaces $\bP_2(x_1, \cdots, x_9)$ carrying an elliptic fibration. The first proposition improves Theorem 2.1 in
case $n(-K_X) = 1,$ compare also (2.11) in [8authors].

\begin{proposition} Let $X$ be a smooth projective 3-fold with $-K_X$ nef. Suppose $n(-K_X) = 1$ and let $f: X \to B = \bP_1$ be
the nef reduction. Then there exists $m_0 \in \bN$ such that $m_0 K_X$ is spanned by global sections and such that the
sections define the map $f.$ In particular $K_X^2 = 0.$ 
\end{proposition}

\proof Let $F$ be the general fiber of $f.$ Then $K_F = K_X \vert F \equiv 0,$ hence there exists $m$ such that $mK_F = \sO_F$
for most fibers. Thus $f_*(-mK_X)$ is a line bundle over $B$ and there is an inclusion
$$ f^*f_*(-mK_X) \to -mK_X.$$
Hence we can write
$$ -mK_X = f^*(\sO_B(a)) + \sum \lambda_i F_i $$
with fiber components $F_i.$ It follows that $ \sum \lambda_i F_i$ is $f$-nef which is only possible if $n \sum \lambda_i F_i =
f^*(\sO_B(b))$ (cut e.g. by a general hyperplane section). Hence we find $m_0$ such that
$$ -m_0 K_X = f^*(\sO_B(c)) $$
and our claim follows. \qed

\begin{theorem} In the situation of (5.1), suppose that there exists a Mori contraction $\varphi: X \to S$ with $\dim S \leq 2.$
Then $-K_S$ is nef. Moreover $$f \times \varphi: X \to B \times S
\simeq \bP_1 \times S$$ is a two-sheeted cover ramified over some $D \in \vert \sO_B(2) \hat \otimes -2K_S \vert.$ The projection $\varphi$
is a conic bundle with discriminant locus $\Delta $ being a member of the linear system $\vert -4K_S \vert$, and the map $f$ is a 
K3-fibration over $B = \bP_1$ with $-K_X = f^*(\sO_B(1)).$ In particular $-K_X$ is spanned by global sections and $-K_X$ is
hermitian semi-positive.
\end{theorem}

Of course, examples as in the theorem exist: just start with $X$ as a smooth two-sheeted covering of $S$,  ramified over $D$ as in the theorem (since $D$
is divisible by 2, the cyclic cover exists). Of course, the existence of a smooth $D \in \vert -K_S \vert$ must be guaranteed.
\\
\\

\proof Let $F$ be a general smooth fiber of $f.$ Since $K_F \equiv 0,$ we must have $\dim \varphi(F) = 2,$ so $S$ is a (smooth) surface, and
$\varphi$ is a conic bundle by Mori's classification. By (2.5) $-K_Y$ is nef.
Let $\Delta$ denote its discriminant locus. Then we have the well-known and easy formula
$$ \varphi_*(K_X^2) \equiv -(4K_S+\Delta).$$
Since $K_X^2 = 0$ by (5.1), $\Delta \equiv -4K_S,$ hence $\Delta = -4K_S$, since $X$ is rationally connected, hence simply connected. 
In particular $\Delta \ne 0$, since $S$ is necessarily a rational surface. \\
Let $l$ be the fiber over a general point of $\Delta.$ Then $l$ is a reducible conic $l = l_1+l_2$ and the $l_i$ are homologous in $X$, since
$\varphi$ is the contraction of an extremal ray. Thus $-K_X \cdot l_i = 1.$ 
Let $d = \deg(f \vert l)$ over a general conic $l = \varphi^{-1}(s),$ then $d \geq 2.$ 
We will show that $d = 2$ so that $-K_X = f^*(\sO_B(1)) $ (use again rational connectedness to pass from numerical equivalence to linear
equivalence). \\
\\
{\bf (I)} We first assume that $K_F = \sO_F.$ \\
Then $f_*(-K_X) $ is a line bundle $\sO_B(a),$ and as in the proof of 5.1, we can write
$$ -K_X = f^*(\sO_B(a)) + \sum_{i=1}^k \lambda_i F_i, $$   
where $F_i$ are fiber components and $\lambda_i$ are positive integers. We take $a$ maximal with such a decomposition and also
note that numerical and linear equivalence coincide, $X$ being simply connected. 
Since $-K_X$ is nef, $\sum \lambda_i F_i = cF$ with a positive rational number $c.$ Since 
$f_*(\sO_X(\sum \lambda_i F_i)) = \sO_B,$ it follows $0 \leq c < 1.$ Thus $a \geq 0,$ again by nefness of $-K_X.$ \\
Suppose next $a = 0.$ 
Since $ \varphi \vert F_i$ is finite, we have $F_i \cdot l_j \geq 1$ (remember that $l_1 $ and $l_2$ are homologous), hence 
$F_i \cdot l \geq 2.$ By virtue of $-K_X \cdot l = 2,$ we must have $k = 1$ and $\lambda_1 = 1,$ so that $-K_X = F_1$. So $h^0(-K_X) = 1$ and the sequence
$$ H^2(\sO_X) \to H^2(\sO_X(F_1)) \to H^2(N_{F_1}) $$
together with $h^2(N_{F_1}) = h^0(K_X \vert F_1)$ shows that $h^2(-K_X) \leq 1.$ Thus $\chi(-K_X) \leq 2.$ On the other hand, Riemann-Roch plus
$K_X^3 = 0$  gives
$$ \chi(-K_X) = 3 \chi(\sO_X) = 3,$$
contradiction. 
Hence $a \geq 1.$ Since $$ 2 = -K_X \cdot l = da + (\sum_{i=1}^{k} \lambda_i F_i) \cdot l, $$ 
we must have $da = 2$ and $k = 0$, therefore $d = 2$ and $  a = 1 $ and
$$ -K_X = f^*(\sO_B(1)).$$ Moreover $f \vert l$ and $\varphi \vert F$ are two-sheeted coverings. Now we consider
$$\tau = f \times \varphi: X \to B \times S,$$
which is a cyclic cover of degree 2, ramified over say $D \subset B \times S.$ Then
$$ K_X = \tau^*(K_{B\times S} +  {{1} \over {2}}D)  $$
together with $K_X = f^*(\sO_S(1))$ proves that $D \in \vert \sO(2) \hat \otimes -2K_S \vert.$ 
Finally, $h =\varphi \vert F \to S$ is a two-sheeted cover which easily shows that $F$ must be K3 (consider $h_*(\sO_F) =
\sO_S \oplus -K_S$ and take cohomology resp.\ take preimages of exceptional curves in $S$).\\
\\ 
{\bf (II)} If $K_F \ne \sO_F,$ then $K_F$ is torsion, write $\lambda K_F = \sO_F.$ Actually $F$ is a hyperelliptic surface or an Enriques surface and 
$\lambda = 2,3,4,6$ in the first resp. $\lambda = 2$ in the second case. \\
\\
{\bf (IIa)} Suppose that $\lambda = 2, $ i.e. $2K_F = \sO_F.$ \\
Arguing as in (I), we have
$$ -2K_X = f^*(\sO_B(a)) + \sum_{i=1}^k \lambda_i F_i.$$ 
As before, $a \geq 0.$ Assuming $a = 0,$ we conclude that 
either $k = 2$ and $-2K_X = F_1+F_2$ or $k = 1$ and $-2K_X = m F_1$   with $m = 1,2.$ 
Now the contradiction is derived in the same way as before: the cohomology groups $H^2(N_{F_1+F_2}) = H^0(K_X \vert F_1+F_2)$ resp.
$H^2(N_{mF_1})$ for $1 \leq k \leq 2$ have dimension at most 2 and 
Riemann-Roch gives $\chi(-2K_X) = 5.$ \\ 
Therefore $a \geq 1$ and the equation
$$ -2K_X = f^*(\sO_B(a)) + \sum_{i=1}^k \lambda_i F_i$$ 
leads to either $da = 2$ and $k = \lambda = 1$ or to $da = 4$ and $k = 0.$ In the first case $a = 1$ and 
$$ -2K_X = f^*(\sO_B(1)) + F_1 = F + F_1,$$
in the second (when $a = 1, d = 4)$
$$ -2 K_X = f^*(\sO_B(1)) = F$$
resp. $$-2 K_X = f^*(\sO_B(2))$$
(when $a = 2, d = 2$). This last case clearly contradicts $K_F \ne \sO_F$.  
In the two remaining cases we have $h^0(-2K_X) = 2$ and $h^2(-2K_X) \leq 2$ thanks to $-2K_X = F+F_1$ resp. $-2K_X = F$ and
$$H^2(\sO_X) \to H^2(\sO_X(F+F_1)) \to H^2(N_{F+F_1}) $$
resp.
$$H^2(\sO_X) \to H^2(\sO_X(F)) \to H^2(N_F).$$
This contradicts
again Riemann-Roch. So $\lambda = 2$ is impossible; in particular $F$ cannot be an Enriques surface.\\
\\
{\bf (IIb)} Now let $\lambda \geq 3,$ in particular $F$ is hyperelliptic. We shall rule out this case. First suppose that $S$
is ruled: $S = \bP(E).$ Then either $S = \bP_1 \times \bP_1$ or $S$ contains a (rational) curve $C$ with $C^2 < 0.$ 
In the first case the covering $F \to S$ produces 2 different fibrations on $F$ which is impossible, in the second $F$ contains a curve $C'$
with $C'^2 < 0,$ which is also absurd. Hence $S = \bP_2,$ in particular the Picard number $\rho(X) = 2.$ Now consider the
relative Albanese map associated with $f,$ 
$$ \sigma: X \rightharpoonup W. $$
$\sigma \vert F$ is the Albanese of the general smooth hyperelliptic surface $F.$ Let $\pi: \hat X \to X$ be a sequence of blow-ups
such that the induced map $\hat \sigma: \hat X \to W$ is holomorphic. Let $A$ be ample on $W$ and consider
$$L = (\pi_*(\sigma^*(A)))^{**}.$$
Then $L$ is a line bundle with $L\vert F $ is nef, non-trivial but not ample.
Since $\rho (X) = 2,$ we can write, at least as $\bQ$-divisors:
$$ L \equiv f^*(\sO_B(a)) \otimes \varphi^*(\sO(b)).$$
Thus $$ L \vert F = \varphi^*(\sO_S(b)),$$ hence $L \vert F$ is ample, trivial or negative, contradiction.

\begin{theorem} In the situation of (5.1) suppose that there is a birational Mori contraction $\varphi: X \to Y.$ Then the following
holds:
\begin{enumerate}
\item $\varphi $ is the blow-up of a smooth curve $C$ in the smooth threefold $Y.$
\item $-K_X = f^*(\sO_B(1))$ and $f$ is a K3-fibration.
\item The normal bundle $N_{C/Y} $ is of the form $N_{C/Y} = L \oplus L$ with some line bundle $L$ on $C.$
\item $-K_Y \cdot C = 2g(C)-2 = \deg L.$
\item If $\deg L > 0,$ then $-K_Y$ is big and nef. 
\item If $\deg L = 0,$ then $-K_Y$ is nef with $K_Y^3 = 0$ and $C$ is elliptic. We have a nef reduction $g: Y \to \bP_2$
such that $-K_Y = g^*(\sO(1))$ and $C$ is a fiber of $g.$
\item If $\deg L < 0,$ then $-K_Y$ is not nef and $C = \bP_1$ with $N_C = \sO(-2) \oplus \sO(-2).$ 
\end{enumerate}
\end{theorem}

\proof The birational map $\varphi$ contracts a prime divisor $E$ either to a point or to a curve. The first alternative however
cannot appear since then $f \vert E$ would have to be finite. So $E$ is contracted to a curve $C$, and $Y$ is automatically smooth
with $\varphi $ being the blow-up of $C$ in $Y.$ \\
Let $l \simeq \bP_1$ be a non-trivial fiber of $\varphi.$ Since $-K_X \cdot l = 1,$ we see with the same methods as in (5.2) that
$$ \deg (f\vert l) = 1 $$
and $$ -K_X = f^*(\sO_B(1)).$$
Denoting $F' = \varphi (F)$ for a fiber $F$ of $f,$ we conclude that $\varphi \vert F: F \to F' $ is an isomorphism. 
Since $H^1(\sO_X) = 0 $ and $H^2(\sO_X(-F)) = H^2(K_X) = H^1(\sO_X) = 0$, the general $F$ is a K3-surface.
Now the exceptional divisor $E$ has two contractions $\varphi \vert E$ and $f \vert E$ so that $E \simeq B \times C.$ In particular we can write
$N_C = L \oplus L.$ \\
Let $$C_b = \{b\} \times C;$$
then $K_X \cdot C_b = 0$ since $C_b$ is contracted by $f.$ From $K_X = \varphi^*(K_Y) + E,$ we deduce
$$ K_Y \cdot C = - E \cdot C_b.$$ 
Now $N_E^* = \sO_{\bP(N^*_C)}(1) = C_b + \varphi_E^*(L^*).$ Hence
$$ -E \cdot C_b =  (C_b + \varphi_E^*(L^*)) \cdot C_b = \deg L^*,$$
and in total 
$$ K_Y \cdot C = \deg L^*. \eqno (1) $$
So the adjunction formula gives
$$ 2g(C)-2 = \deg L. \eqno (2) $$
From the exact sequences
$$ 0 \la N_{C_b/E} = \sO \la N_{C_b/X} \la N_{E/X} \vert C_b = L \la 0$$
and
$$ 0 \la N_{C_b/F} \la N_{C_b/X} \la N_{F/X} \vert C_b = \sO \la 0$$ we obtain
$$ N_{C_b/F} = L. \eqno (3)$$
We also notice 
$$ \sO_X(F) = \varphi^*(\sO_Y(F')) - E $$
and
$$ N_{F/X} =  N_{F'/Y} - \sO_{F'}(C), $$
thus
$$ N_{F'/Y} = \sO_{F'}(C). \eqno (4) $$ 
Since $-K_Y$ is nef on every curve $\ne C,$ the bundle $-K_Y$ is nef precisely when $\deg L \geq 0 $ by virtue of (1).
The equation $$0 = K_X^3 = (\varphi^*(K_Y) + E)^3 $$ 
gives 
$$ K_Y^3= \deg L^*. \eqno (5)$$
Finally we observe that because of $-K_X = F,$ the formula
$$ -K_Y = F' \eqno (6)$$
holds.\\
\\
{\bf Case I:} $\deg L > 0.$ \\
Then $-K_Y$ is nef and big by (5) and the previous remark. The normal bundle $N_{F'/Y} $ is big and nef by (6).\\
\\
{\bf Case II:} $\deg L = 0.$ \\
Here $-K_Y$ is nef with $K_Y^3 = 0.$ The normal bundle $N_{F'/Y} $ is effective (and  nef); actually $N_{F'/Y} = \sO_{F'}(C).$ 
Furthermore it is clear that $n(-K_Y) = 2.$ In fact, $E \cap F$ is an elliptic curve $l$  for general $F$ and $0 = K_X \cdot l =
K_Y \cdot \varphi(l) $ giving a 2-dimensional $K_Y$-trivial family of elliptic curves on $Y.$ To be more precise, we consider
the exact sequence
$$ 0 \to H^0(\sI_C \otimes (- K_Y)) \to H^0(-K_Y) \to H^0(-K_Y \vert C) \to H^1(I_C \otimes (- K_Y)).$$
Now $I_C \otimes (- K_Y) = \varphi_*(-K_X),$ hence $h^0(I_C \otimes (- K_Y)) = 2$ and $h^1(I_C \otimes (- K_Y)) = 0,$ as one checks immediately
(in fact, $h^1(-K_X) = 0)$.
Using the normal bundle sequence for $C \subset F' \subset Y$ and $N_{C/F'} = \sO,$ the normal bundle $N_{C/Y}$ must be trivial or the non-split
extension by two trivial bundles, 
hence $-K_Y \vert C = \sO_C.$ Putting this into the exact sequence, we conclude that $h^0(-K_Y) = 3$ and $-K_Y$ is spanned.
Let $\tilde h: Y \to \bP_2$ be the associated map (which contracts $C$)
and let $h: Y \to T$ be its Stein factorisation. Then $-K_Y^2 = \tilde h^{-1}(x)$
for any $x \in \bP_2,$ on the other hand $K_Y^2 = C.$ We conclude that then $h$ must be an isomorphism, hence $T = \bP_2.$ \\ 
\\
{\bf Case III:} $\deg L < 0.$ \\
Then (2) shows that $C = \bP_1$ and that $N_C = \sO(-2) \oplus \sO(-2).$ In this case $-K_Y$ is not nef. \\
 
\begin{example} {\rm (1)  Let $Z = \bP(E)$ be a $\bP_3$-bundle over $\bP_1$ such that $Z$ is Fano; let $\pi$ denote its projection. 
We suppose furthermore that $-K_Z \otimes \pi^*(\sO(-1)$ is generated by global sections (the existence of a smooth section
would be sufficient). Take a general
$$ X \in \vert -K_Z \otimes \pi^*(\sO(-1) \vert.$$ 
Let $\psi: Z \to W$ denote the second projection and let $f = \pi \vert X;$ $\varphi = \psi \vert X.$ 
Then $$ -K_X = f^*(\sO(1))$$
the general fiber of $f$ is a quartic in $\bP_3$ and $n(-K_X) = 1.$ The condition on $-K_Z \otimes \pi^*(\sO(-1)$ is translated into
$$ S^4E \otimes \det E^* \otimes \sO(-1) $$
being generated by sections. If we write
$$ E = \bigoplus \sO(a_i) $$
with $a_1 \geq \ldots \geq a_4,$ then this condition comes down to 
$$ 4a_4 \geq \sum a_i + 1 = c_1(E) + 1.\eqno (*)$$ 
This implies that up to a twist $E$ can only be of type $(0,0,0,0)$ or of type $(1,0,0,0).$
In  the first case $Z = \bP_1 \times \bP_3.$ 
It is immediately checked that $\varphi: X \to \bP_3$ is birational and that $\varphi$ is the blow-up of a curve $C$ with $\deg C = 16 $
and $g(C) = 33.$ \\
In the second case, the second contraction $\psi: Z \to W = \bP_4$ is the blow-up in a plane $S$ and $X$ is the blow-up of a 
quartic $W'$ in $\bP_4$ along $W' \cap S.$ \\ 
\\
(2) In order to get an example for 5.3(6), we take a threefold $Y$ such that $-K_Y$ is spanned by global sections and such that $K_Y^3 = 0.$ 
Then let $X \subset \bP_1 \times Y$ be a general element of $\vert \sO(1) \hat \otimes -K_Y \vert.$ \\
\\
(3) At the moment we do not have an example for 5.3(7).}

\end{example}

\begin{theorem} Rationally connected threefolds $X$ with the following properties are bounded modulo boundedness for threefolds $Y$
with big and nef $-K_Y$ resp.\ threefolds with $-K_Y$ nef and $n(-K_Y) = 3.$  
\begin{itemize}
\item $-K_X$ is nef;
\item $n(-K_X) = 1;$ 
\item $X$ does not admit a contraction which is of type $(-2,-2).$  
\end{itemize}
\end{theorem}

The boundedness - actually classification - of threefolds with $-K$ big and nef is under investigation at the moment; contractions
of type $(-2,-2)$ are excluded just for technical reasons and reasons of length. We come back to that in a separate paper.

\proof (1) First suppose that $X$ a contraction $\varphi: X \to S$ with  $\dim S \leq 2.$ Then Theorem 5.2 applies and, using the
notations of (5.2), it only remains to bound the bundle $E,$ i.e. to bound $S.$ But since $K_X^2 = 0,$ we actually have 
$-K_S$ nef (see proof of 2.2)), thus $S$ is bounded. \\
(2) If $\varphi$ is birational to the threefold $Y$, then (5.3) applies. Using the notations of (5.3) and ruling out the
case of a $(-2,-2)$-contraction, we have either have $(-K_Y)^3 > 0 $ or $K_Y^3  = 0$ with $C$ elliptic. In the second case, we have
an elliptic fibration $g: Y \to \bP_2$ and $C$ is a fiber. In order to proceed by induction on the Picard number, we want to apply
Theorem 6.9. For that we verify that $Y$ does not admit a $(-2,-2)$-contraction. In fact, if $D$ is the exceptional divisor
of such a contraction, then $D$ meets only singular fibers of $g$ and therefore $D \cap C = \emptyset. $ But then $D$ 
defines already a $(-2,-2)$-contraction on $X$ which was ruled out by assumption.  \\
If however $(-K_Y)^3 > 0,$ then $Y$ is bounded (by assumption), hence it remains to bound $C$ for fixed $Y.$ Let $a = K_Y^3.$
From
$$ 0 = K_X^3 = K_Y^3 + 3 \varphi^*K_Y \cdot E^2 + E^3 $$
and $\varphi^*(K_Y) \cdot E^2 = - K_Y \cdot C$ and $E^3 = - c_1(N_{C/Y})$ we obtain
$$ 2 K_Y \cdot C + 2g-2 = K_Y^3 = a. \eqno (*)$$
Here $g$ denotes the genus of $C.$ Now consider $\vert -mK_Y \vert $ for some large $m$ and the associated birational embedding
$$ \psi: Y \to Y' \subset \bP_N.$$
Let $\lambda: Y \to Y'$ be the birational part of $\psi$ and put $C' = \lambda (C).$ Also notice $K_Y = \lambda^*(K_{Y'}).$ 
By the theory of Hilbert schemes it suffices to bound the degree of $C'$ and the genus of $C'$ ( $ = g)$ (if $\dim \lambda (C) = 0,$
then $(*)$ proves that $C = \bP_1$; on the other hand, its normal bundle is ample due to (5.3), so this case cannot occur).
Due to $(*)$ we only need to bound $-K_Y \cdot C.$ But this follows from $K_X^2 = 0,$ i.e. $K_Y^2 = C,$ hence $-K_Y \cdot C = (-K_Y)^3.$

\qed 

\section{Rationally connected threefolds II}

\begin{setup} {\rm  We are now turning to the case of rationally connected threefolds $X$ with $-K_X$ nef and $n(-K_X) = 2.$ 
Here we have a holomorphic elliptic fibration $f: X \to B$ to a normal projective
surface $B$ by (2.1).  Since $-mK_X = f^*(G)$ for some ample line bundle $G$, there exists 
an effective $\bQ$-divisor $D$ on $B$ such that 
$(B,D)$ is log-terminal and $K_X  \equiv f^*(K_B+D)$ [Na88,0.4]. In particular $B$ has only quotient singularities. 
Again we consider a Mori contraction $\varphi: X \to Y.$ \\
We note that by Riemann-Roch $\chi(-K_X) = 3$ since $K_X^3 = 0.$ Since $K_X^2 \ne 0$ then by Kodaira vanishing $H^2(-K_X) = 0,$
therefore
$$ h^0(-K_X) \geq 3. \eqno (6.1.1) $$}
\end{setup}

\begin{theorem} Suppose in (6.1) that $\dim Y = 1.$ Then $B = \bP_2$ and the general fiber $F$ of $\varphi$ is $\bP_1 \times \bP_1$
or a del Pezzo surface. In particular $X$ is a ramified cover over $\bP_1 \times \bP_2$ of degree $d$ at most 8. Moreover we have
$$ -K_X = f^*(\sO_{\bP_2}(a)) $$
with $1 \leq a \leq 2$ and $a = 2$ can only happen when $F = \bP_1 \times \bP_1.$ The elliptic fibration $f$ is equidimensional.
Finally, if $a = 1,$ then $K_F^2  = d$ and $8 \geq K_F^2 \geq 2.$  
\end{theorem} 

\proof Notice that $f \vert F$ is finite and $\varphi$ is finite on
every $f$-fiber $F_f.$ Therefore the equidimensionality of $f$ is clear: if $F_0$ is a 2-dimensional fiber of $f,$ then $F_0 \cap F$ would be
a curve.
\\
(1) First we show that $\varphi$ cannot be a $\bP_2$-bundle. This already settles the second  assertion.
Suppose $X = \bP(E)$ with a rank 3-bundle $E$ over $Y = \bP_1.$ Then $K_X^3 = 0$ translates into
$$ c_1(S^3E \otimes \det E^* \otimes \sO(2)) = 0$$
which is absurd. \\
(2) Next observe that $\rho (B) = 1$ and also the group of Weil divisors modulo linear equivalence is $\zed$; just because $\rho(X) = 2.$  
Let $C \subset F$ be a smooth rational curve with $C^2 = 0$ and let $C' = h(C),$ where
$h = f \vert F.$ Then $C'$ is a moving rational curve in $B$ not meeting the singularities of $B.$ Hence $H^0(\omega_B) = 0$,
so that $H^2(\sO_B) = 0.$ Since anyway $H^1(\sO_B) = 0,$ the surface $B$ has only rational singularities (and in particular again
$K_B$ is $\bQ$-Cartier). This is also clear from the fact that $B$ has only log-terminal singularities.  
\\
(3) Consider the torsion free sheaf $f_*(-K_X)$ and let $\sL$ be its reflexive hull. Then at least outside a set of codimension at least 2 in $X$ we can 
write 
$$ -K_X = f^*(\sL) + D_0$$
where $f_*(\sO_X(D_0)) = \sO_B.$ This can also be considered as an equation of $\bQ$-Cartier divisors on all of $X$, resp. as
$-K_X = f^*(\sL)^{**} + D_0$ on all of $X.$  
Now we have
$$2 = -K_X \cdot C = (\deg h_C) (\sL \cdot C') + D_0 \cdot C.$$
Thus we are in one of the following cases:
\begin{enumerate}
\item $D_0 \cdot C =  \deg h_C = \sL \cdot C' = 1$;
\item $D_0 \cdot C = 0; \ \deg h_C = 2; \ \sL \cdot C' = 1$;
\item $D_0 \cdot C = 0; \ \deg h_C = 1; \ \sL \cdot C' = 2$.
\end{enumerate}
But we know that $K_X \equiv f^*(K_X+D),$ hence as a $\bQ$-divisor, $D_0 \equiv f^*(G).$ This shows that $D_0 \cdot C > 0$ unless
$D_0 = 0.$ Hence in cases (2) and (3) we have $D_0 = 0.$    
By (6.1.1) we have $h^0(-K_X) \geq 3,$
hence $h^0(\sL) \geq 3$. \\
\noindent \vskip .2cm
Now suppose that we are in one of the first two cases.
Therefore always $\sL \cdot C' = 1.$ Hence any effective Weil 
divisor on $B$ is a positive (integer) multiple of $\sL$, having in mind that $C'$ lies in the regular part of $B.$ Considering the
exact sequence
$$ H^0(\sL \otimes \sO_B(-C')) \to H^0(\sL) \to H^0(\sL \vert C') $$
and having in mind that $C' = \bP_1$ (the $C'$ are linearly equivalent, thus the general $C'$ is smooth), we conclude that 
$\sL = \sO_B(C')$ and $h^0(\sL) = 3.$ In particular $\sL$ is locally free, $C'^2 = 1,$ $K_B = \sO_B(-3C')$ and thus
$B = \bP_2$ (use e.g. Fujita's $\Delta$-genus, [Fu90];[BS95]). \\
\vskip .2cm \noindent
(4) We still have to consider the case $\sL \cdot C' = 2.$ 
\vskip .2cm \noindent
(4a) Arguing in the same way suppose first that $H^0(\sL \otimes \sO_B(-C')) \ne 0.$
Then we either have $\sL = \sO_B(C')$ with $C'^2 = 2$
or $\sL = \sO_B(2C')$ with $C'^2 = 1.$ In the first case $K_B$ is divisible by 2 and $B$ is a quadric cone; in the second $K_B$ is
divisible by 3 and $B = \bP_2.$ We rule out the case that $B$ is a quadric cone as follows. We can write 
$$ h_*(\omega_{F/B}) = \sE \oplus \sO_B  \; ;$$
here $\sE$ is a reflexive sheaf at least (it is not clear whether $h$ is flat; we will not care about that). 
To proceed we first notice that
$$ H^1(h_*(\omega_{F/B}) \otimes \sO_B(-2)) \ne 0. \eqno (*)$$
On the other hand, we are going to prove
$$ H^1(h_*(\omega_{F/B}) \otimes \sO(-2)) = 0. \eqno (**)$$
This comes down to the vanishing
$$ H^1(F,\omega_{F/B} \otimes h^*(\sO_B(-2))) = 0.$$
Now $-K_X = f^*(\sO_B(1))$ in our situation, hence $\omega_{F/B} = h^*(\sO_B(1)).$ Since $\omega_F = h^*(\sO_B(-1)),$
our claim therefore comes down to prove that
$$H^1(F,\omega_F) = 0.$$
This is however clear by duality and $(**)$ is verified so that in case(4a) $B$ must be a plane. \\
\\
(4b) If now $H^0(\sL \otimes \sO_B(-C')) = 0,$ then $h^0(\sL) = h^0(\sL \vert C') = 3.$ So $H^0(\sL)$ defines a meromorphic map
$$ g: B \rightharpoonup \bP_2$$ which is holomorphic near $C'$; moreover $g(C')$ is a conic and thus $C'^2 = 4.$ \\ 
Suppose that $\sO_B(C')$ generates ${\rm Pic}(B) = \zed.$ Then, having in mind that $\sL \cdot C' = 2, $ we must have
$$ \sL = \sO_B({{1} \over {2}}C').$$
Hence $\sO_B(C') = g^*(\sO_{\bP_2}(2)) $ near $C'$ and since $C'^2 = 4,$ we conclude that $g$ must be generically $1:1.$ Thus $g$ is an 
isomorphism outside
the finite set of indeterminacies. But now the linear system $\vert C' \vert $ defines an embedding into $\bP_5$ which factors via $g$ over
the Veronese embedding of $\bP_2$ and therefore proves that $g$ is an isomorphism. 
\\
If $\sO_B(C')$ is not the ample generator $\sO_B(1)$, then necessarily $\sO_B(C') = \sO_B(2),$ and then $\sL = \sO_B(1)$, in particular
$\sL$ is locally free. Now $c_1(\sL)^2 = 1$ and hence again $g$ is an isomorphism.\\
\vskip .2cm \noindent
(5) We now show that $D_0 = 0$ which has to be proved only in the first case.
Supposing $D_0 \ne 0,$ we have
$$ -K_X = f^*(\sO_B(1)) + D_0$$
and we first consider the case that $F$ is as del Pezzo surface (different from the quadric).
Then take a $(-1)$-curve $l \subset F$, and we obtain
$$ 1 = -K_X \cdot l = f^*(\sO_B(1)) \cdot l + D_0 \cdot l. $$
Since $l$ and $C$ are homologous in $X$, we have $D_0 \cdot l > 0,$ contradiction. \\
If $F = \bP_1 \times \bP_1,$ then using the equation of $\bQ$-divisors $D_0 = f^*(\sO_B(b))$ and $D_0 \vert F = (1,1)$ 
(since $D_0  \cdot C = 1$),
we obtain $D_0 \vert F = h^*(\sO_B(1)) $, so that $b = 1.$ This is absurd.\\
\noindent \vskip .2cm
(6) Finally, if $a = 2$ then $-K_F$ is divisible by 2 which is only possible if $F = \bP_1 \times \bP_1$. The last statement of the theorem is
clear.

\qed

\begin{example} {\rm (1) Let $g: X \to \bP_1 \times \bP_2$ be the cyclic cover of degree 2, ramified over a smooth divisor of type
$(4,2)$. Then $-K_X = p_2^*(\sO(2)),$ so that $-K_X$ is nef but not big. Moreover $p_2 $ is the nef reduction and $p_1$ is a
quadric bundle. \\
(2) We modify the previous example by taking $R$ of type $(4,4)$. Then $-K_X = p_2^*(\sO_B(1))$ and $p_1$ is a del Pezzo fibration
whose general fiber $F$ has $K_F^2 = 2$ (hence $F$ is $\bP_2$ blown up in 7 points). \\
(3) Let $h: X \to \bP_1 \times \bP_2$ be the cyclic cover of degree 3 ramified over a smooth divisor of type (3,3). Then
$-K_X = p_2^*(\sO_B(1))$ and $p_1$ is a del Pezzo fibration with $K_F^2 = 3.$ \\
In all examples it is easily checked that indeed $b_2(X) = 2$ so that $p_1$ is the contraction of an extremal ray. }
\end{example}

\begin{theorem} In the setup (6.1) suppose that $\dim Y = 2.$ Let $\Delta$ be the discriminant locus of the conic bundle $\varphi: X
\to Y$. Suppose that $\Delta \ne 0.$ Then 
\begin{enumerate}
\item The nef reduction (given by $\vert -mK_X \vert$ for suitable large $m$) is equidimensional;
\item $-(4K_Y+\Delta) $ is nef; and
\item either $f \times \varphi: X \to B \times Y$ is an embedding; and $B = \bP_2$ 
\item or $f \times \varphi: X \to B \times Y$ is a 2:1-covering onto its image and $B = \bP_2.$
\end{enumerate}
\end{theorem}

\proof By (2.1) $-mK_X$ is spanned for suitable large $m$ and therefore defines the nef reduction. $l = l_y, y \in Y$ will always denote a 
fiber of $\varphi$ and we set $l' = l_y' = f(l_y).$ 
\noindent \vskip .2cm
(A) Assume that there is a 2-dimensional fiber component $S$ of $f$.
Since $h^0(-K_X) > 0,$ we can write
$$ -K_X = aS + D_0 + E$$ with effective divisors $D_0$ and $E$ such that $E = f^*(E')$ and $f_*(\sO_X(D_0)) = \sO_B$  and with a positive integer $a.$ 
In fact, choose a point $p \in S$ and a section $s \in H^0(-K_X)$ vanishing at $p;$ then notice that $K_X \vert S \equiv 0$ so that $s \vert S = 0.$
Let $a$ be the vanishing order and then consider $-K_X - aS.$ \\ 
Now let $l$ be a general fiber of $\varphi.$ 
Since $\Delta \ne 0$, we must have $S \cdot l \geq 2.$ Thus $a = 1,$ $S \cdot l = 2,$  $E \cdot l = D_0 \cdot l = 0$ so that $D_0 = \varphi^*(\tilde D)$
and $E = \varphi^*(\tilde E).$ 
Let $F$ be a general 
fiber of $f.$ Then $(E + D_0) \cdot F = 0, $ hence $(\tilde E + \tilde D) \cdot \varphi(F) = 0.$ But $\kappa ( \tilde E + \tilde D) = 2$ since $\kappa (-K_X) = 2$,
contradiction to
the fact that $\varphi(F)$ moves in $Y.$ This proves (1).\\
Notice that as a consequence of (1) the general $l'$ does not meet any singularity of $B.$ 

\noindent \vskip .2cm \noindent
(B) Claim (2) follows from 
$$ -(4K_Y+\Delta) \cdot C = K_X^2 \cdot \varphi^{-1}(C) \geq 0.$$ 
Note that $-(4K_Y+\Delta) \ne 0,$ since $K_X^2 \ne 0.$
\noindent \vskip .2cm \noindent
(C) Approaching (3) and (4) we write as in (A)
$$ -K_X = E + D_0 =
f^*(E') + D_0.$$ The second equation is an equation of $\bQ$-divisors; in terms of sheaves it reads
$$ \sO_X(-K_X) = f^*(\sO_B(E'))^{**} \otimes \sO_X(D_0).$$
Notice that $E' \ne 0;$ actually $h^0(\sO_B(E')) \geq 3.$ 
By intersecting with irreducible components of reduced conics
we obtain $D_0 \cdot l = 0$ for all $l.$ 
Since $E \cdot l = 2,$ we must have $\deg (f \vert l) \leq 2.$ \\
Next we show that $\rho(B) \leq 2$ and actually $\rho(B) = 1$ if $(l')_{y\in Y}$ is a 2-dimensional family. 
In fact, if $(l')$ is 1-dimensional, then take a general curve $C \subset Y$ and set $X_C = \varphi^{-1}(C).$ Clearly $X_C$ 
projects onto $B.$ But all fiber components of $X_C$ are homologous in $X$ (not in $X_C$), thus $\rho(B) \leq 2.$ 
If $(l')$ is 2-dimensional, then choose $x_0 \in B$ general and pick an irreducible curve $C \subset Y$ such that $x_0 \in l'_y$
for all $y \in C.$ Then $X_C \to B$ contracts some curve and therefore $\rho (B) = 1.$ \\

\vskip .2cm \noindent 
(D) First we assume $\rho (B) = 1.$ 
Let $\sO_B(1)$ be the ample generator on $B.$ 
We write
$$ \sO_B(l') = \sO_B(a) \, , \, \, \sO_B(E') = \sO_B(b) \ {\rm and}  \ -K_B = \sO_B(c) $$
with a positive integer $a$  and positive rational numbers $b,c.$ 
Let $d = c_1(\sO_B(1))^2. $\\
\vskip .2cm \noindent
(D.1) Suppose that $E' \cdot l' = 2.$ 
Suppose first that $b < a.$ Then $E' \cdot l' = 2$ implies $h^0(\sO_B(E')) \leq 3,$ hence $h^0(\sO_B(E')) = h^0(-K_X) = 3.$ 
To see this, we are going to show that for $E'$ and $l'$ general, $E' \cap l'$ is contained in the smooth locus of $l'.$ Once we know this,
$h^0(\sO_B(E') \vert l') \leq 3$ by considering the normalization of $l'.$ If $\vert E' \vert $ has no fixed components, the intersection
statement is clear. Otherwise we could write $E' = M + E''$ with $M$ fixed. Since $B$ is easily seen to be $\bQ$-factorial (since $\rho (B) 
= 1$) and since
$M \cdot l' > 0,$ we conclude that $\sO_B(M) = \sO_B(E'')$ which is absurd. \\
To continue, let 
$$g: B \rightharpoonup \bP_2$$ 
be the associated rational map. Then $g$ is holomorphic near the general $l'$ since sections on $l'$ lift to $B.$ 
Moreover $g$ maps $l'$ biholomorphic onto a conic in $\bP_2,$ in particular $l'$ is smooth.
Since two conics in $\bP_2$ meet in 4 points generally, we must 
have $l'^2 \leq 4,$ in other words
$$ a^2 d \leq 4.$$ Thus $a = 1$ and $d \leq 4$ or $a = 2$ and $d = 1.$ 
Now the equation $E' \cdot l' = 2$ translates into 
$$ 2 = abd.$$ Putting things together, either $a = 1$ or $(a,b) = (2,1)$ and if $a = 1,$ then $d = 4$ and $b = {{1} \over {2}}.$ 
In this case we compute the $\Delta$-genus
$$ \Delta(\sO_B(1)) = 2 + d - h^0(\sO_B(1)) = 6 - h^0(\sO_B(1)).$$
Namely, Kawamata-Viehweg and Riemann-Roch give 
$$ h^0(\sO_B(1)) = \chi(\sO_B(1)) = 6.$$
Hence $\Delta(\sO_B(1)) = 0$ and Fujita'classification (see e.g. [BS95,3.1.2]) implies that there is a smooth rational curve
$C \subset B$ such that $B $ is $\bP(\sO_C \oplus \sO_B(1) \vert C)$ with the zero-section blown down. Then an easy explicit
calculation on $\bP(\sO_C \oplus \sO_B(1) \vert C)$ shows that the invariant $e = 4$. Let $\pi: \hat B \to B$ be the canonical desingularization, i.e. 
$\hat B = \bP(\sO \oplus \sO(-4))$ over $\bP_1.$ Let $C_0 $ be the negative section and $F$ a fiber of $\pi.$ Then 
$\pi^*(E') = \alpha C_0 + 2F$ for some $\alpha$ (pull-back as reflexive sheaf) since $\pi^*(\sO_B(1)) = C_0 + 4F.$ 
We also know that $f \vert l$ is biholomorphic, $f(l) $ being a smooth rational curve. Thus $f \times \varphi$ is an embedding and
$X \cap (B \times y)$ is defined by $\sO_B(E').$ In other words, $l' \in \vert \sO_B(E') \vert.$ But there is no irreducible reduced
member in $\vert E' \vert,$ which follows immediately by considering $\vert \alpha C_0 + 2F \vert$ for any $\alpha$. 
\vskip .2cm \noindent
Now to the case that $a = 2$ and $b = 1.$ Again we compute 
$$ \Delta(\sO_B(1)) = 0.$$
Arguing as before and using again $h^0(\sO_B(1)) = 3,$ we see that $e = 1$ and therefore $B = \bP_2.$\\

If $b \geq a,$ then $a = 1$ and $(b,d) = (1,2)$ or $(2,1)$. Suppose we know that $l'$ is smooth for general $l'.$ Then  
the adjunction formula for $l' \subset B$ gives
$$ -2 = -cad + a^2d = d(1-c).$$ 
Thus in case $(b,d) = (1,2)$ we get $c = 2$ and $B$ is the quadric cone  by the (generalized) Kobayashi-Ochiai theorem. If $(b,d) = (2,1),$ then
$c = 3$ and $B = \bP_2$ by the same quotation. We exclude the case of the quadric cone: it is clear that $f \times \varphi$
is an embedding and that $X \subset B \times Y$ is defined by $\sO_B(1) \hat \otimes -K_Y.$ In particular $X$ is Cartier in 
$B \times Y$. But $X$ is smooth and must meet the singular locus of $B \times Y,$ contradiction. \\
It remains to check the smoothness of $l'$. If $c > a = 1,$ then 
$$ 0 = H^1(\sO_B) \to H^1(\sO_{l'}) \to H^2(\sO_B(-l')) = H^0(\sO_B(l') + K_B) = 0$$
proves $H^1(\sO_{l'}) = 0$ so that $l' = \bP_1.$ 
If $c \leq 1,$ then Riemann-Roch for $\chi(\sO_B(1))$ yields $c = 1.$ Thus $B$ is Gorenstein (with ample anticanonical class). 
Then it is well-known that $-K_B$ is spanned (e.g. by classification). Hence $\sO_B(l')$ is spanned and $l'$ is smooth (of course
this immediately contradicts $l' \in \vert -K_B \vert$). \\
\vskip .2cm \noindent 
(D.2) If $E' \cdot l' = 1,$ then from $h^0(\sO_B(l')) \geq 3,$ we obtain $h^0(\sO_B(E'-l')) \ne 0.$ Unless $E' = l'$  we can write $E' = l' + R$
which clearly contradicts $\rho (B) = 1.$ Now we have $E' = l'$ which means $b = a.$ Moreover $E'\cdot l' = 1$ translates into
$ 1 = abd,$ hence $a = b = d = 1.$ Computing $\chi(\sO_B(1))$ we get $c = 1$ or $3$, the case $c = 1$ being excluded as in (D.1). Hence $c = 3$
and $B = \bP_2.$

\vskip .2cm \noindent
(E) Finally we treat the case $\rho (B) = 2.$ Then the family $(l')$ is 1-dimensional. 

\vskip .2cm \noindent
(E.1) Case $\deg f \vert l = 1$: \\
Let $l'$ be general and choose $b \in l'$
such that $f^{-1}(b)$ is elliptic. Let
$$C = \varphi(f^{-1}(l')).$$
Notice that $\dim C = 1$ by our assumption that the family $(l')$ is 1-dimensional. By the degree assumption, $C$ is a
(possibly singular) elliptic curve. Also, every fiber of $f$ over $l'$ must be elliptic since it is mapped onto $C.$ 
Let $$X_C = \varphi^{-1}(C).$$
Then both $\varphi({\rm Sing}(X_C))$ and $f({\rm Sing}(X_C))$ are finite. Hence ${\rm Sing}(X_C) $ is finite and thus $X_C$
is normal. Thus $C$ and $l'$ are smooth and $X_C \simeq C \times l'.$ Therefore we are in one of the the two following situations. 
\noindent \vskip .2cm
(I) $f: X \to B$ is an elliptic fiber bundle outside a finite set in $B.$ 
\noindent \vskip .2cm
(II) There are finitely many $l_1', \ldots , l_k'$ such that all $f$-fibers over every $l'_i$ are singular and
moreover every $l'$ different from $l'_1, \ldots , l'_k$ is disjoint from $\bigcup l_j'.$ 
\noindent \vskip .2cm
In case (I) we have $R^1f_*(\sO_X) = \sO_B$ in codimension 1 whence $q(X) > 0$ as $f$ is equidimensional
contradicting the rational connectedness of $X.$ \\
In case (II) we consider the graph of the family $(l')$ and deduce immediately the existence of a map $g: B \to T = \bP_1$ contracting
all $l'.$ Since all fibers of $\varphi$ are contracted by $X \to T,$ there is a map $h: Y \to T$ such that 
$$ g \circ f = h \circ \varphi.$$ Since $B$ is $\bQ$-factorial and
$\rho (B) = 2,$ all fibers of $g$ must be irreducible. Also, no fiber can be multiple, otherwise
by base change $\varphi$ would have multiple fibers in codimension 1. Thus all fibers of $g$ are irreducible reduced and therefore
$\Delta = \emptyset, $ contradiction.

\vskip .2cm \noindent 
(E.2) Case  $\deg (f \vert l) = 2$: \\
Then $E' \cdot l' = 1$ and
the $l'$ form a 1-dimensional family. Let $l'$ be general, so that $l'$ is contained in the regular part
of $B.$ 
Consider $X_l = f^*(l').$ Then $\dim \varphi (X_l) = 1$ and 
$N_{l/X_l} = \sO_l.$ Now consider the exact sequence 
$$ 0 \to f^*(N^*_{l'/B}) \to N^*_{l/X} \to N^*_{l/X_l} \to 0.$$
Since $N^*_{l/X} = \sO_l \oplus \sO_l,$ we conclude that $N^*_{l'/B} = \sO_{l'},$ in particular $(l')^2 = 0.$
Since $-K_B-D$ is ample, $K_B \cdot l' < 0.$ Hence we conclude $l' \simeq \bP_1$. Thus $H^0(\sO_B(l'))$ defines a
holomorphic map $g: B \to C \simeq \bP_1$ contracting $l'.$ 
The general fiber $F$ of the induced map $X \to \bP_1 = C$ has $-K$ nef (and not ample) and is therefore either $\bP_1 \times {\rm elliptic} $ or 
$\bP_2$ blown up in 9 points. In any case $g $ induces a map $h: Y \to C $ such that $h \circ \varphi = g \circ f.$ 
Now $g$ has irreducible fibers due to $\rho (B) = 2$ and also $g$ cannot have multiple fibers as in (E.1 (II)). Hence $g$ is a $\bP_1$-bundle. 
If we write $B = \bP(V) \to C$, then $X = \bP(h^*(V)) \to Y,$ and thus
$\varphi$ is not a proper conic bundle. 

\vskip .2cm \noindent (F) Finally we observe that in case $\deg f \times \varphi = 1$, this map is an embedding by Zariski's Main Theorem.
In case of degree 2, we already saw that $B = \bP_2$. Also it follows that $f \vert l$ is a degree $2$ covering for all smooth $l$, resp. 
an isomorphism on all components of singular conics. Therefore $f \times \varphi$ is a {\it  covering} of degree 2.   

\qed

\begin{example} {\rm (1) Let $X \subset \bP_2 \times \bP_2$ be a smooth divisor of type $(3,2)$. Then $-K_X$ is nef, and one projection
defines an elliptic fibration while the other is a conic bundle.  It is not difficult to see directly that the
discriminant locus $\Delta$ is always non-empty.\\
(2) If we take in (1) $X$ of type $(3,1)$, then instead of a conic bundle we have a $ \bP_1$-bundle coming
from a vector bundle, a situation we study next.
Of course we can take $X$ more generally as a smooth hypersurface in $Y \times \bP_2$ of type $(-K_Y,1)$, where the requirement on $Y$ is just that
$\vert  -K_Y \hat \otimes \sO(2) \vert$ contains a smooth member. \\
(3) To get an example with a $2:1$-covering, let $Z = \bP(T_{\bP_2})$ with projection $p: Z \to \bP_2$  (i.e. $Z \subset \bP_2 \times \bP_2$ has degree
$(1,1)$.
Let $G = \sO_Z(1) \otimes p^*(\sO(1));$ then $G$ defines the second contraction $q: Z \to \bP_2$ of the
Fano manifold $Z \subset \bP_2 \times \bP_2.$ Now take a smooth element $R \in \vert 4G \otimes p^*(\sO(2)) \vert$. Let $h: X \to Z$ be the 2:1-covering ramified over $R.$ Then $K_X = h^*p^*(\sO(-1)),$ so that 
$-K_X$ is nef with an elliptic fibration $X \to \bP_2.$ The map $q \circ h$ defines a conic bundle 
structure over $\bP_2.$ }
\end{example}

\begin{theorem} In (6.1) suppose that $\dim Y = 2$ and that $\varphi: X \to Y$ has discriminant locus $\Delta = \emptyset. $
Then $\varphi$ is a $\bP_1$-bundle and of the form $X = \bP(E)$ with a rank 2 vector bundle $E$ over $Y$. 
In particular $-K_Y$ is nef. Furthermore:
\begin{enumerate}
\item $B = \bP_2, \bP_1 \times \bP_1$ or $\bP_2$ blown up in one point.
\item If $B = \bP_2$, then $X \subset \bP_2 \times Y$ is given by $\sO(1) \hat \otimes (-K_Y)$ and $K_Y^2 > 0$ 
\item If $B = \bP_1 \times \bP_1$ or $\bP_2$ blown up in one point, then 
$Y$ is $\bP_2$ blown up in 9 points in such a way that $Y$ carries an elliptic fibration 
$g: Y \to D = \bP_1$ and such that
$E = g^*(\sO(a) \oplus \sO)$ with $a = 0,1.$ 
\end{enumerate}
\end{theorem}

\proof Since $Y$ is a smooth rational surface, $H^2(Y,\sO_Y^*) = H^3(Y,\zed)$ is torsion free; hence every analytic $\bP_1$-bundle 
over $Y$ is of the form $\bP(V).$ Then $K_X^3 = 0$ translates into
$$ 3K_Y^2 = 4c_2(V) - c_1^2(V).  \eqno (*)$$
\vskip .2cm \noindent
(A) First we prove that $f$ is equidimensional. So suppose to the contrary that $S$ is a 2-dimensional fiber component. 
Then $S$ must be a section of $\varphi$: consider a general curve $C \subset Y$ and let $X_C = \varphi^{-1}(C).$ Then $f \vert X_C$ 
is generically finite; on the other hand $S$ projects onto $Y$, so $S \cap X_C$ is a curve, i.e. $X_C$ contains a contractible curve,
which is necessarily a section. Hence $S$ is a section of $\varphi$. (Of course $S \to Y$ is finite!) \\
This section corresponds to an exact sequence (after a suitable twist of $V$)
$$ 0 \to \sO_Y \to V \to L \to 0 \eqno (S_1) $$
such that $S = \bP(L).$ Notice that $h^0(L) = 0,$ since $S$ does not move, and $\sO_X(S) = \sO_{\bP(E)}(1).$ 
Since $f$ contracts the surface $S$, we have $K_X \cdot S = 0,$ in particular $$K_X \cdot S \cdot \varphi^*(H) = 0$$
for all ample divisors $H$ on $Y.$ This comes down to 
$$L \cdot H - K_Y \cdot H = 0,$$
hence $L = K_Y.$ 
Putting this into $(*)$, it follows $3 L^2 = - L^2,$ hence $L^2 = 0 = K_Y^2.$ In particular $Y$ is $\bP_2$ blown up in 9 points in almost general 
position. \\
Suppose first that $(S_1)$ splits. Since $-mK_X$ is spanned for suitable $m,$ we easily see that $-mK_Y$ is spanned and therefore
$Y$ has an elliptic fibration. This elliptic fibration induces the elliptic fibration on $X$ and it is clear that $f$
cannot have a 2-dimensional fiber. It remains to rule out the case that $Y$ does not carry an elliptic fibration and 
at the same time $(S_1)$ does not split. Then however, taking symmetric powers of $(S_1),$ $h^0(S^m(V \otimes -K_Y))$ can grow
at most linearly, contradicting the spannedness of $-mK_X$ and $K_X^2 \ne 0.$ 
Thus $f$ does not have a 2-dimensional fiber.    
\vskip .2cm \noindent
We consider the ruling family $(l_y)_{y  \in Y}$ in $X.$ As in (6.4) the image family $( l'_y) = (f(l_y))$ in $B$ is either 2- or 1-dimensional.  
Using the decomposition $-K_X = E + D_0$ as in (6.4) and the obvious fact that
$E \cdot l_y > 0,$ the degree $\deg (f \vert l_y)$ is 1 or 2. Also we have $\rho (B) \leq 2$ and $\rho(B) = 1$ if $(l')$ is 
2-dimensional. 

\vskip .2cm \noindent
(B) Case $\rho (B) = 1$: \\
Then we can argue exactly as in (6.4) and conclude that
$B = \bP_2$. Moreover $f \times \varphi$ is an embedding or a degree 2 covering over its image in $B \times Y.$ \\
First suppose that $f \times \varphi$ is an embedding. Then
$X \in \vert \sO(1) \hat \otimes H \vert$. Here $H$ is some line bundle on $Y$. 
$X$ cannot be of degree $2$
in $\bP_2$ because then $\varphi$ would be a conic bundle with $\Delta = \emptyset;$ on the other hand the reducible conics in $\bP_2$ 
have codimension 1 in the parameter space $\bP_5$ of all conics, so we must have $X = l \times Y$ with a fixed conic $l$. Then however
$Y$ must have an elliptic fibration $g: Y \to C$ and $B = l \times C$ contradicting $\rho (B) = 1.$ \\
So $X$ is linearly embedded in the
trivial $\bP_2$-bundle over $Y.$ Therefore after a suitable twist $V$ is a quotient of the trivial 
rank 3-bundle, so that there is an exact sequence
$$  0 \to M \to \sO_Y^3 \to V \to 0 \eqno (S') $$
with a line bundle $M$ on $Y$. Notice $M = \det V^*.$ 

Computing $K_X$ by the adjunction formula for $X \subset B \times Y$ and
via the $\bP_1$-bundle structure shows by comparison that $H = \det V.$ 
Then $(S')$ gives $c_1^2(V) = c_2(V).$ 
Let $l_b = \varphi(f^{-1}(b)),$ an elliptic curve for general $b \in B.$ We conclude that
$l_b \in \vert \det V \vert.$ This also shows $h^0(V) = 3$ thanks to $h^1(\det V^*) = 0.$ 
Now the adjunction formula yields
$$ 0 = K_Y \cdot l_b + c_1^2(V),$$
hence $c_1^2(V) = - K_Y \cdot \det V.$ 
On the other hand, Riemann-Roch gives 
$$ \chi (V) = -{{1}\over {2}} c_1^2(V) - {{1} \over {2}} c_1(V) \cdot K_Y + 2 = 2,$$
and $\chi (V) = h^0(V) - h^1(V).$ Thus $h^1(V) = 1$ and therefore $h^2(\det V^*) = h^0(\det V + K_Y) = 1.$ 
Take $0 \ne D \in \vert K_Y + \det V \vert.$ \\ 
First we assume $K_Y^2 > 0.$ Then from $D \cdot (-K_Y) = 0$ we 
deduce
$$ D = \sum a_i C_i$$
with $a_i \geq 0$ and $C_i$ some $(-2)$-curves. Therefore
$$ \det V = -K_Y + \sum a_i C_i.$$
Since $\det V $ is nef, this is only possible if all $a_i = 0.$ Hence $\det V = -K_Y$  if $K_Y^2 > 0.$ \\
Now suppose $K_Y^2 = 0.$
So $c_1(V)^2 = c_2(V) = K_Y^2 = 0,$ and $Y$ has an elliptic fibration $g: Y \to C;$ moreover $\det V = - a K_Y.$ 
Using (S') and taking into account possible multiple fibers we find a rank 2-bundle $V'$ over $C$ and a line bundle
$L$ over $Y$ such that $V = g^*(V') \otimes L.$ In particular we conclude that $B = \bP(V'),$ a contradiction to our assumption that
$\rho (B) = 1.$ \\
If $f \times \varphi$ has degree $2,$ then consider its image $X' \subset B \times Y.$ Consider the ramification
divisor $R.$ Then $R \cdot p_Y^{-1}(y) = 2,$ i.e.\ $R \to Y$ is a degree 2 covering and it must be ramified since
$Y$ is simply connected. Over the ramification points in $Y$ we therefore have reducible conics. Contradiction.

\noindent \vskip .2cm \noindent
(C) Case $\rho (B) = 2$:  \\
Thus $(l')$ is 1-dimensional. Arguing as in (E.1) of the proof of (6.4), we obtain a map $g: B \to C = \bP_1$ and a map $h: Y \to C$ such that 
$ g \circ f = h \circ \varphi.$ 
Now we argue as in (B) above to get $B = \bP(V').$ The nefness of $-K_X$ is equivalent to the nefness of
$$ g^*(V' \otimes {{\det V'^*} \over {2}}) \otimes {{-K_B } \over {2}}.$$
This means that --up to normalization-- $V' = \sO \oplus \sO $ resp. $V' = \sO \oplus \sO(1)$ so that $B = \bP_1 \times \bP_1$
resp. $\bP_2(x).$  
 
\qed

\begin{proposition} In (6.1) suppose that $\dim Y = 3.$ Let $E$ denote the exceptional divisor and suppose that $\dim \varphi (E) = 0.$ 
Then $-K_Y$ is big and nef and we are in one of the following situations.
\begin{enumerate}
\item $E = \bP_2$ with normal bundle $N_E = \sO(-1);$ $Y$ is smooth with $(-K_Y)^3 = 8;$ $B = \bP_2$ with $\deg f_E = 1 $ or $ 4.$ 
\item $E = \bP_2$ with normal bundle $N_E = \sO(-2);$ $(-K_Y)^3 = {{1} \over {2}};$ $B = \bP_2$ with $f_E$ an isomorphism, i.e. $E$ is a section of $f.$ 
\item $E = \bP_1 \times \bP_1$ with $N_E = \sO(-1,-1);$ $B = \bP_1 \times \bP_1$ or $\bP_2$ and $f_E$ is an isomorphism resp. $\deg f_E = 2.$
\item $E = Q_0$, the quadric cone, with $N_E = \sO(-1);$ $B = Q_0$ or $\bP_2$ with $f_E$ an isomorphism resp. $\deg f_E = 2.$  
\end{enumerate}
\end{proposition} 

\proof 
It is clear that $-K_Y$ is nef and the classification of $(E,N_E)$ is provided by [Mo82]; from this information also the computation
of $(-K_Y)^3$ is clear. It remains to determine $B$ and $\deg f_E$. 
It is clear that $E$ maps onto $B$ and also that $f_E$ is finite. In particular any 2-dimensional
fiber component of $f$ is disjoint from $E.$  
Again we write in codimension 1:
$$ -K_X = f^*(\sL) + D \eqno (*)$$
where $D$ is the contribution from the multiple fibers. 
\vskip .2cm \noindent 
(1) Here $-K_X \vert E = \sO(2).$ 
Restricting to a general line
$l \subset E$ we obtain either $\deg f \vert l = 1$, $\sL \cdot l' = 1$ and $D \cdot l = 1$ or $\deg f \vert l = 1$, $\sL \cdot l' = 2$, $D = 0$
or $\deg f  \vert l = 2$, $\sL \cdot l' = 1$ and $D = 0.$ 
Here $l' = f(l).$ 
\\ Suppose first that $\deg f \vert l = 1$. By Bertini $f^{-1}(l')$ must be irreducible for general $l'$, hence $\deg f_E = 1$ and $B = \bP_2.$
Now $f$ cannot have multiple fibers and therefore $D = 0,$ ruling out the first case.
\\ So $\deg f \vert l = 2$ and $L \cdot l' = 1.$
Then $\sO_B(l') \in {\rm Pic}(B) = \zed$ and the $\bQ$-line bundle $\sL$ can be written as
$\sL = \sO_B(\mu l')$ with a positive rational number $\mu.$ Since we know by (6.1.1) that $h^0(\sL) \geq 3,$ the exact sequence
$$ H^0(B,\sO_B(-l')) \to H^0(B,\sL) \to H^0(\sL \vert l')$$
together with $\sL \cdot l' = 1$ shows that $\mu \geq 1.$ Then $1 = \sL \cdot l' = \mu l'^2,$ hence $\mu = 1.$ So $\sL = \sO_B(l')$
and the three sections of $\sL$ give an isomorphism $B \to \bP_2.$ Squaring $(*)$ we see that $\deg f_E = 4.$ 
\vskip .2cm \noindent(2) Here $-K_X \vert E = \sO(1),$ so that $D = 0$ and $\deg f \vert l = 1.$ Again we conclude $\deg f_E = 1$, in particular $B = \bP_2.$ 
\vskip .2cm \noindent (3) Since $-K_X \vert E = \sO(1,1)$ and $h^0(L) \geq 3,$ we have $D = 0.$  Restricting to a ruling line $l$, we obtain $\deg f \vert l = 1$
and $L \cdot l' = 1.$ Distinguishing the cases $\rho (B) = 1$ and $2$ and using both ruling families, similar calculations as in (1) lead
to $B = \bP_2$ resp. $\bP_1 \times \bP_1.$ In the first case $\deg f_E = 2,$ in the second $f_E$ is an isomorphism. 
\vskip .2cm \noindent (4) The case of a quadric cone is similar. 
\qed

\begin{proposition} In (6.1) suppose that $\dim Y = 3;$ let $E$ be the exceptional divisor and
suppose that $\dim \varphi (E) = 1;$ i.e. $Y$ is smooth and $\varphi$ is the blow-up of a smooth curve
$C.$ 
Then we are in one of the following 2 cases:
\begin{enumerate}
\item $\dim f(E) \leq 1$.
Then $-K_Y$ is nef with $n(-K_Y) = 2$ unless $\varphi$ is of type $(-2,-2);$ 
if in case $-K_Y$ nef $g: Y \to B'$ denotes the nef
reduction of $Y,$ then there exists a birational map $\tau: B \to B'$ such that $\tau \circ f =
g \circ \varphi $ and $C$ is a smooth (elliptic) fiber of $g.$   
\item $\dim f(E) = 2$. Then $B$ is a rational ruled surface and there exists another contraction
$\psi: X \to Y'$ such that either $\dim Y \leq 2$ or $\dim Y = 3$ but the exceptional divisor $E$ 
does not project onto $B.$ 
\end{enumerate} 
\end{proposition} 

\proof (1) Suppose first that $\dim f(E) \leq 1.$ 
Let $A' = f(E)$ and let $A$ be the image of the
Stein factorization of $E \to A';$ in particular $A$ is a smooth rational curve. 
Then $E$ admits two different projections, hence $E \simeq C \times A$ (at least after finite \'etale cover which we
will ignore for simplicity of notations).
We write $E = \bP(N^*),$
where $N^*$ is the conormal bundle of $C \subset Y.$ Then $N$ decomposes: $N = L \oplus L$ with
a line bundle $L$ on $C.$ Let $g$ be the genus of $C$ and let $C_0 = C \times a $ for some $a \in A.$
Then by the standard theory of ruled surfaces
and adjunction we can write
$$ -K_X \vert E = C_0 + (2-2g + \deg L) F,$$
with $F$ a fiber of $\varphi \vert E.$ 
Since $-K_X$ is nef, we conclude that $\deg L \geq 2g-2.$ 
Now $C_0$ is contracted by $f$ so that $-K_X \cdot C_0 = 0.$ This is translated into 
$\deg L = 2g-2.$ Since $C_0$ is contained in a fiber of $f$, $L$ cannot be ample and thus 
we have $g \leq 1.$ If $g = 0$
then $L = \sO(-2).$ This is the $(-2,-2)$ case and $-K_Y$ is nef except on $ C.$
Thus we suppose that $g = 1.$ Then $K_Y \cdot C = 0$, so that $-K_Y$ is nef and $K_Y^3 = 0 $. Clearly $n(-K_Y) = 2;$ 
let $g: Y \to B'$ be the nef reduction. Now $f \vert E $ is just the projection onto $A = A' = \bP_1$,
$-K_X$ being nef. Moreover $N_{A/B} $ is negative and $B'$ is just the blow-down of $A \subset B.$ Thus $g$ contracts $C$ and
$C$ is just a fiber of $g.$ 
\\ (2) Let $l$ be a general fiber of $\varphi \vert E.$ 
Similarly as in the proof of (6.4)(1) we see that $E$ is disjoint from all potential 2-dimensional fiber components. Now use the decomposition $-K_X = f^*(\sL) + D$ in codimension 1 and intersect with $l$
to obtain $D = 0$, $\deg f \vert l = 1$ and $\sL \cdot l' = 1$ for $l' = f(l).$  
\\ If $\rho (B) = 1,$ then $l'^2 > 0$ so that $f_E^*(l') $ is ample. But $l \subset f^*_E(l')$, thus Bertini gives a contradiction.  
\\ So $\rho (B) = 2.$ Then $l'^2 = 0$ and the linear system $\vert l' \vert $ gives a morphism $\tau: B \to T$ to a smooth curve $T$ 
such that $\tau \circ f_E = h \circ \varphi_E $ with some covering $h: C \to T.$ Since $\tau$ is flat with all fibers irreducible and reduced, 
it must be a $\bP_1$-bundle. So $B$ is ruled and of course rational. Let $p: X \to \bP_1$ be some projection, $p = q \circ f_E$ with  
$q: B \to \bP_1$ a $\bP_1$-bundle structure. Then $K_X$ is not $q$-nef, hence there exists a relative contraction $\psi: X \to Y'$. 
This is the contraction we are looking for. 
\qed

\begin{theorem} Smooth rationally connected threefolds subject to the following conditions are bounded modulo boundedness
of threefolds $Y$ with $-K_Y$ nef and $n(-K_Y) = 3:$
\begin{itemize}
\item $-K_X$ is nef
\item $X$ has a contraction not of type $(-2,-2)$
\item $n(-K_X) = 2.$ 
\end{itemize}
\end{theorem}

\proof Let $\varphi: X \to Y$ be a contraction. If $\dim Y = 1,$ we have boundedness by (6.2). In case $\dim Y = 2$ and the discriminant
locus $\Delta \ne \emptyset,$ 
(6.4) clearly gives boundedness once we can bound $Y.$ But $-(4K_Y + \Delta) $ is nef, and therefore $Y$ is bounded. 
If $\Delta = \emptyset, $ then boundedness follows from (6.6). 
So we are reduced to the case that $\varphi$ is birational. By assumption we may assume that $\varphi$ is not of type $(-2,-2).$
If $\varphi$ is a $(-1,-2)$-contraction, then by [DPS93, 3.5/3.6/3.7] $X$ is canonically isomorphic to a $\bQ$-Fano threefold 
with a fixed type of singularities, hence $X$ is bounded. So we may assume that $-K_Y$ is nef. If the exceptional divisor $E$ contracts
to a point, then (6.8) applies and we conclude by assumption. Finally let $\dim \varphi(E) = 1.$ If $\dim f(E) = 2,$ then
by (6.8) we can switch to another contraction of a different type and work with that one. If $\dim f(E) \leq 1$ and if $\varphi$ is not
of type $(-2,-2)$, then we conclude by induction on $\rho (X).$
\qed

\section{Rationally connected threefolds III}  

\begin{setup}  {\rm The last case to deal with is $n(-K_X) = 3$ for a rationally connected threefold $X$ with $-K_X$ nef. Hence there is 
no covering family $(C_t)$
such that $K_X \cdot C_t = 0.$ Of course this is the case when $(-K_X)^3 > 0,$ i.e. $-K_X$ is big and nef. Then $-mK_X$ is generated for suitable
large $m$ and the curves $C$ with $-K_X \cdot C = 0$ are just those which are contracted by the morphism associated with $\vert -mK_X \vert.$
So we shall assume $K_X^3 = 0.$ From the proof of theorem 2.1 we know that in this case $K_X^2 \ne 0$  i.~e.\ $\nu(-K_X)=2$.
}
\end{setup}

\subsection{The structure of the anticanonical system}

\begin{proposition} \label{prop72}
Let $X$ be a smooth projective rationally connected threefold with  
$-K_X$ nef, $n(-K_X)=3$ and $\nu(-K_X)=2$. Then 
the anticanonical system $|- K_X|$ has non-empty fixed part 
$A$. It induces a fibration $f: X \to \bP_1$. If $F$ is a
fiber of $f$ then $|-K_X| = A + |kF|$ with $k \geq 2$.
Furthermore $A^3 = A^2 \cdot F = 0$. 
\end{proposition}

\proof We know $h^0(-K_X) \geq 3$ by Riemann-Roch and Kawamata-Viehweg 
vanishing.
Assume that $|-K_X|$ is without fixed part. Then for two general members
$D,D' \in |-K_X|$ we have $D \cdot D' = C$ with an effective curve
$C$ on $D$ and $-K_X \cdot C = 0$ as $K_X^3 = 0$. This implies also that
$-K_X \cdot C' = 0$ for any component $C'$ of $C$. Now we assumed $n(-K_X)=3$
and therefore the $-K$-trival curves do not cover $X$. Hence, as $D$ moves, we 
conclude that no component of $C$ moves on $D$.\\
An easy calculation using $q(X)=0$ shows that $h^0(\sO_D(C)) \geq 2$
which implies $h^0(\sO_D) \geq 2$ and $D$ has at least two 
connected components $D=P+Q$. Let $H \subset X$ be very ample and 
let $D_H, P_H, Q_H$ denote the restrictions to $H$. First we note that
$P_H$ and $Q_H$ are nef divisors on $H$: If $C$ is a curve in $H$ which
is contained in $P_H$ then $Q_H \cdot C = 0$ as $P_H$ and $Q_H$ do not 
meet and therefore $P_H \cdot C = D_H \cdot C \geq 0$. Now $D_H^2 > 0$
as $K_X^2 \neq 0$ and we may assume that $P_H^2 >0$. As $Q_H$ is
orthogonal to $P_H$ and nef the Hodge index theorem implies that 
$P_H$ and $Q_H$ are proportional. In particular $P_H^2 = 0$ as 
$P_H \cdot Q_H = 0$ which gives the desired contradiction.\\
Now write $|-K_X| = A + |B|$ with non-empty fixed part
$A$ and movable part $B$ and $h^0(B) = h^0(-K_X) \geq 3$.\\
As $-K_X$ is nef we know that $K_X^2 \cdot A \geq 0$ and $K_X^2 \cdot B \geq
0$. Now $0 = -K_X^3 = K_X^2 \cdot (A + B)$ gives $K_X^2 \cdot A = 
K_X^2 \cdot B = 0$. From this we further conclude that $-K_X \cdot
(A \cdot B + B^2)=0$. As B moves $A \cdot B$ and $B^2$ are 
effective cycles which implies $-K_X \cdot A \cdot B = -K_X \cdot B^2
=0$. From the last equation and $-K_X^2 \cdot A = 0$ we finally deduce
that $-K_X \cdot A^2=0$.\\
In the next step $-K_X \cdot B^2 = 0$ in combination with $n(-K_X)=3$ gives 
the further structure of $B$. As $B$ has
obviously no fixed part we may repeat the argumentation in the first part of
the proof with $B$ playing the role of $D$ to prove that $B$ has at least
two connected components. Now $B_H$ has no fixed part and is therefore nef.
So following the proof above a few lines more shows that $B^2 \neq 0$ 
leads to a contradiction. Hence $B^2 = B_H^2 =0$ and as $|B_H|$
has at most finitely many base points, some multiple of $B_H$ defines a
morphism sending $H$ to a curve. In particular every connected component
of $B_H$ (being nef) is equivalent to a rational multiple of a general
fiber of this map so that two connected components of $B$ are equivalent
up to a rational factor. As $B$ has at least two such components this implies
that some multiple of $B$ is generated by its global sections.\\
Let $f: X \to \bP_1$ be the Stein factorization of the morphism defined by
$|mB|$ for sufficiently large $m$ and let $F$ be a fiber. Then $B$ is 
equivalent to a rational multiple of $F$. As $|B|$ has no fixed part we
can actually conclude that $B$ consists of $k$ fibers with $k$ an integer
and as $B$ has at least two connected components $k \geq 2$.
Finally $-K_X \cdot A^2 = 0$ gives $A^3 + A^2 \cdot B = 0$ and $K_X^3=0$ 
gives $A^3 + 3 A^2 \cdot B=0$ (as $B^2=0$) which together imply that
$A^2 \cdot B = A^3 = 0$.
\qed

\begin{corollary}
In this situation $\rho(X) \geq 3$.
\end{corollary}

\proof Assume that $\rho(X) = 2$ and let $H$ be an ample divisor. As 
$A$ and $F$ are not proportional we can write $H = \alpha A + \beta F$
with some rational numbers $\alpha, \beta$. As $A^3 = A^2 \cdot F = F^2 = 0$
we calculate $H^3 = 0$ which is impossible.
\qed 

Next we study the fibration defined by $-K_X$.

\begin{lemma}
Let $F$ be a general fiber of $f$. Then $F$ is 
a smooth surface with $-K_F$ nef and effective, $K_F^2=0$ and $n(-K_F)=2$.
\end{lemma}

\proof The canonical bundle formula shows that $-K_F = -K_{X|F} = A_{|F}$. As 
$A^2 \cdot F = 0$ we have $K_F^2=0$. Furthermore $n(-K_F)=2$ otherwise
$F$ and therefore $X$ would be covered by $K_X$-trivial curves.
\qed

Now we consider the different cases corresponding to the type of a
general fiber $F$. There are three different cases: 
\begin{enumerate}
\item $F$ is $\bP_2$ blown
up in 9 (not necessarily distinct) points without elliptic fibration and the 
anticanonical divisor is a smooth elliptic curve. 
\item As before, but this time the anticanonical divisor
consists of rational curves. 
\item $F$ is a $\bP_1$-bundle over an elliptic curve (with $n(-K_F) = 2)$.
\end{enumerate}

\begin{proposition}
In the setting of proposition \ref{prop72} write $|-K_X| = A_1 + A' + |kF|$ ($k\geq2$) 
where $A_1, A'$ are effective divisors with ${A_1}_{|F} = -K_F$ resp.\ $A'_{|F} = 0$
for a general fiber $F$. Furthermore, let $\sL = \sO_X(-A'-kF) = \sO_X(K_X+A_1)$. Then:
\begin{enumerate}
\item $h^0(\sO_{A_1})= 1 + h^1(R^1 f_* \sL)$ and $h^2(\sO_{A_1})=0$
\item $h^1(\sO_{A_1})= h^0(\sO_{\bP_1}(l+k-2) + h^0(R^1 f_* \sL) \geq 1$ where
the nonnegative number $l$ is defined via $f_* (\sO(-A') = \sO_{\bP_1}(-l)$. In particular, 
$h^1(\sO_{A_1})= 1$ iff $k=2$, $A'=0$ and $h^0(R^1 f_* \sL)=0$.
\end{enumerate}
\end{proposition}

\begin{remark} If $F$ is rational, $R^1 f_* \sL$ is torsion. So in this case $h^0(R^1 f_* \sL)=0$
iff $R^1 f_* \sL=0$.
\end{remark}

\proof Everthing follows from direct computations: As $X$ is rationally connected $H^i(\sO_X)=0$ for 
$i>0$. The sequence
$$ 0 \to \sO_X(-A_1) \to \sO_X \to \sO_{A_1} \to 0 $$
then implies that
$$ h^0(\sO_{A_1}) = 1 + h^1(O_X(-A_1))$$
$$ h^1(\sO_{A_1}) = h^2(O_X(-A_1))$$
$$ h^2(\sO_{A_1}) = h^3(O_X(-A_1))$$
By duality, $h^i(O_X(-A_1)) = h^{3-i}(\sO_X(-A'-kF)) = h^{3-i}(\sL)$, in particular
$h^2(\sO_{A_1}) = h^0(-A'-kF) = 0$. The sheaf $f_* \sL$ is torsion free hence
locally free and of rank one as $\sL_{|F}=\sO_F$. In fact $f_* \sL = \sO_{\bP_1}(-l-k)$
where $l$ is the number of fibers containing components of $A'$. For the further calculations
we use the Leray spectral sequence which collapses at the $E_2$ level as the base is a curve. So
$$ h^1(\sO_{A_1}) = h^1(X,\sL) = h^0(R^1 f_* \sL) + h^1 (f_* \sL)$$
and $h^1(f_* \sL) = h^1(\sO_{\bP_1}(-l-k) = h^0(\sO_{\bP_1}(l+k-2)$ which gives the claim. Finally
$$ h^0(\sO_{A_1}) = 1 + h^2(X,\sL) = 1 + h^0(R^2 f_* \sL) + h^1 (R^1 f_* \sL)$$
and $R^2 f_* \sL = 0$ as $h^2(\sL_{|F}) = h^2(\sO_F) = h^0(-K_F) = 0$ for every fiber $F$, again
by duality.
\qed

{\it For the rest of the paper we will concentrate on the case where a general fiber
$F$ is $\bP_2$ blown up in 9 points such that $|-K_F| = C$ a smooth elliptic curve.}

As $-K_X$ is nef and not numerically trivial we have $\kappa(X) = -\infty$.
In particular $K_X$ is not nef and we can study $X$ using some Mori contraction
$\varphi$. If $\varphi: X \to X'$ is birational then $X'$ is again smooth and we
stay within the small list of \cite{Mo82a}; in particular we do not 
encounter small contractions or flips. We also keep the fact that the
structure of $|-K|$ is very special. The final outcome is some Mori fiber
space $\varphi: X \to Y$.

\subsection{The case where $F$ is rational}

\subsubsection{The setup}

\begin{proposition} \label{prop75}
Consider as above a smooth projective rationally connected threefold $X$ with $-K_X$ nef,
$n(-K_X)=3$, $\nu(-K_X)=2$ and let $f:X \to \bP_1$ be the fibration induced
by $|-K_X| = A + |kF|$, $k \geq 2$ and $F$ a general fiber. We further assume
that $F$ is $\bP_2$ blown up in 9 points such that $-K_F$ is nef and
$|-K_F| = C$ a smooth elliptic curve. Then $k=2$, $A = C \times \bP_1$ 
and $f$ restricted to $A$ is the second projection. 
\end{proposition}

\proof As $-K_F = A_{|F} = C$ a irreducible reduced curve we find a
divisor $A_1$ which occurs in $A$ with multiplicity one and the rest $A'$
does not meet $F$. Furthermore the restriction of $f$ to $A_1$ is an 
elliptic fibration and the anticanonical bundle $-K_{A_1} = (A'+kF)_{|A_1}
\geq kF_{|A_1}$ contains $f_{|A_1}^* \sO(2)$ which will imply our
other assertions:\\
Let $\nu: \hat{A}_1\to A_1$ be the normalization and let $\mu: \tilde{A}_1
\to \hat{A}_1$ be the minimal desingularization. Let $h:\bar{A}_1 \to \bP_1$ be
a relative minimal model of the induced elliptic fibration $g:\tilde{A}_1
\to \bP_1$ i.~e.\ we take the successive blow-down $\lambda$ of ($-1$)-curves 
contained in fibers. Computing the (anti-)canonical bundles we get
$$ -K_{\hat{A}_1} = \nu^* (-K_{A_1}) + Z$$
with some effective Weil divisor $Z$ supported on the zero locus of the
conductor ideal and
$$ -K_{\tilde{A}_1} = \mu^* (-K_{\hat{A}_1}) + E_1$$
with some effective divisor $E_1$ and finally
$$ -K_{\tilde{A}_1} + E_2 = \lambda^* (-K_{\bar{A}_1})$$
with another effective divisor $E_2$. In particular we still have that 
$-K_{\bar{A}_1} \geq h^* \sO(2)$ and the dual of the relative dualizing
sheaf for $h$ is effective. By \cite[Theorem III.18.2]{BPV84} this is only 
possible
if all smooth fibers are isomorphic and the only singular fibers 
are multiple fibers of the form $m_i F_i$. Then the weak canonical bundle
formula for elliptic fibrations \cite[Corollary V.12.3]{BPV84} shows that
$$ K_{\bar{A}_1}=h^* L + \sum (m_i -1)F_i $$
with some line bundle $L$ of degree
$$\deg L = \chi(\sO_{\bar{A}_1}) - 2 \chi(\sO_{\bP_1}) = 1 - q(\bar{A}_1) - 2$$
which reads in total as
$$ -K_{\bar{A}_1} = h^* \sO(1+q(\bar{A}_1)) + \sum (1 - m_i) F_i $$
By the argument above and the Leray spectral sequence we conclude that 
$q(\bar{A}_1) \leq 1$. Together with $-K_{\bar{A}_1} \geq h^* \sO(2)$ 
this shows that there are no multiple fibers and $h$ is a $C$-bundle. But 
the base is $\bP_1$ hence $\bar{A}_1 = C \times \bP_1$. In particular 
$-K_{\bar{A}_1} = h^*\sO(2)$ and every inequality above was in fact an 
equality. This implies that
$$ \bar{A}_1 = \tilde{A}_1 = \hat{A}_1 = A_1 = C \times \bP_1$$
-- here we use that $C \times \bP_1$ has no curves with negative self intersection;
also $Z=0$ implies that $A_1$ is regular in codimension 1. As we already know that
$A_1$ is Cohen-Macaulay, it is normal. By Proposition 7.5 now $q(A_1)=1$ implies 
$A' = 0$ and $k = 2$ which settles the claim.
\qed

\begin{corollary} Under the assumptions of \ref{prop75} we have 
$|-mK_X| = mA + |2m F|$.
\end{corollary}

\proof This follows immediately from the fact that $A$ has just one
component, $A_{|F} = -K_F$ and $\kappa(-K_F) = 0$.
\qed

\subsubsection{Running the MMP -- birational contractions}

We start with divisorial Mori contractions $\varphi: X \to X'$ where
$X$ has the following special structure provided by proposition 7.7: 

\begin{definition} Let $X$ be a smooth projective threefold. We say that 
$X$ has structure \emph{(A)} if there exists a fibration $f: X \to \bP_1$
such that a general fiber $F$ is a smooth rational surface with $-K_F$ nef,
$|-K_F|$ contains a smooth elliptic curve $C$ and $-K_X = A + 2F$
where $A \cong C \times \bP_1$ and $f_{|A}$ is the second projection.
Furthermore we require that $-K_X$ is almost nef which means that there are
at most finitely many rational curves $D_i$ such that $-K_X \cdot D \geq 0$
for all curves $D \neq D_i$.
\end{definition}

\begin{remark}
Let $\varphi: X \to X'$ be a divisorial contraction. By \cite[Prop.\ 
3.3]{DPS93}
if $-K_X$ is nef then in general $-K_{X'}$ is merely almost nef (and
this conclusion also holds if we just require $-K_X$ to be almost nef).
\end{remark}

\begin{proposition}
Assume that $X$ has structure (A). Then there is no birational Mori 
contraction 
$\varphi:X \to X'$ where a divisor is mapped to a point.
\end{proposition}

\proof Assume that we have such a contraction and let $E$ be the exceptional 
divisor. We
first treat the case $E \cong \bP_2$. By Mori's list we have two possibilities 
for the 
normal bundle: either $\sO_E(E) = \sO_E(-1)$ or $\sO_E(-2)$. As $E$ itself is 
not fibered
we must have $F \cdot E = 0$ and $E$ is contained in some fiber $F_0$ of $f$. 
Obviously $E$
is different from $A$, therefore $A \cdot E$ is an effective cycle of curves 
and must
be either the elliptic curve $C_0 = A \cdot F_0$ or zero. On the other hand 
computing 
the canonical bundle $K_E = (K_X + E)_{|E} = (-A + E)_{|E}$ we see that 
$A_{|E}$ has 
degree 1 or 2. Contradiction.\\
Another case is $E \cong Q$ a singular quadric in $\bP_3$. By Mori's 
classification the
normal bundle is $\sO_E(-1)$ in this case and we can conclude as above.\\
The last case is $E \cong \bP_1 \times \bP_1$ with normal bundle 
$\sO_E(-1,-1)$ hence
$-K_{X|E} = \sO_E(1,1)$. If $l$ is a (general) fiber of the first projection 
of $E$ we have 
$-K_X \cdot l = 1$. As $E \neq A$ we know that $l\not\subset A$ and therefore 
$A \cdot l \geq 0$. This gives the only possibility $F \cdot l = 0$, $A \cdot 
l = 1$
which implies that $F_{|E} = \sO_E(a,0)$ with $a \geq 0$ hence $A_{|E} = 
\sO_E(1-2a,1)$.
But this is an effective cycle on $E$ which means that $1-2a \geq 0$ therefore 
$a=0$ and
$F_{|E}$ is trivial. So $E$ is contained in some special fiber $F_0$ and 
$A_{|E}$ is either the
elliptic curve $C_0 = A_{|F_0}$ or zero. Both cases are not possible as 
$A_{|E} = \sO_E(1,1)$.
\qed

\begin{proposition} \label{prop710}
Assume that $X$ has structure (A) and let $\varphi: X \to X'$ be a birational 
Mori contraction
where the exceptional divisor $E$ is mapped to a curve. Then one of the
following two cases occurs:
\begin{description}
\item (i) $E \neq A$ and we have $-K_{X'} = A' + 2F'$ with $A' = \varphi(A)
\cong A, F' = \varphi(F)$ and either $F' \cong F$ or $F \to F'$ is the
blow-down of some (--1)-curves in $F$. The fibration $f$ factors as
$f = f' \circ \varphi$ and $f': X' \to \bP_1$ gives $X'$ the structure (A).
\item (ii) $E = A$ and we have $-K_{X'} = 2G$ with $G = \varphi(F) \cong F$ 
and $\varphi$ is the blow-up of the elliptic curve $G^2$ which is in $|-K_G|$. 
In this case $X'$ has structure (O).
\end{description}
\end{proposition}

\begin{definition} Let $X$ be a smooth projective threefold. We say that 
$X$ has structure \emph{(O)} if $-K_X = 2G$ is almost nef where $G$ is a 
smooth rational surface with $-K_G$ nef, $|-K_G|$ contains a smooth elliptic 
curve $D = G^2$ and $G$ moves in a linear system without fixed components.
\end{definition}

\begin{re}
In case (ii) of proposition \ref{prop710} $f$ does not factor via $\varphi$ as the 
images of two fibers $G_1, G_2$ meet in $D = \varphi(E)$. Since $A \cdot l = -1$ 
this case can only occur if $A$ is not nef.
\end{re}

\proof The extremal ray corresponding to the contraction is generated by a
rational curve $l$ with $-K_X \cdot l = 1$. In fact $\varphi$ is the
blow-up of the smooth curve $D \subset X'$ and $l$ is a fiber of 
the $\bP_1$-bundle $\pi = \varphi_{|E}$. We also know that $X'$ is smooth.\\
We first treat the case $E \neq A$. Then $A \cdot l \geq 0$ and the 
intersection
numbers are $A \cdot l = 1$ and $F \cdot l = 0$. Hence if we pick some
$l_0$ it must be contained in a fiber $F_0$. We now have two different 
subcases:
$\dim f_{|E} = 0$ or $1$. If $f_{|E}$ maps onto $\bP_1$ then $E \cdot F$ is 
non-empty and in fact $F_{|E}=b\:\!l$ with $b>0$. Restricted to a general 
fiber 
$\varphi_{|F}$ blows down some (--1)-curves $l$ in $F$ each of them meeting
$C = -K_F$ transversally in one point. In particular all the curves $l$ meet
$A$ transversally in one point which implies that $\varphi_{|A}$ is an
isomorphism.\\
If $f(E)$ is a point, $E$ is contained in some special fiber $F_0$ and in
particular $A \cdot E$ is either $C_0 = -K_{F_0}$ or zero and in fact it
must be $C_0$ as $A \cdot l = 1$. This also shows
that $\varphi(A) \cong A$. The rest is obvious because $\varphi$ is an
isomorphism outside $F_0$.\\
The other case is $E = A$. Since $A = C \times \bP_1$ we know that $l$ is a 
fiber of the first projection of $A$ and $F \cdot l = 1$ as $A_{|F}=C$.
As we contract the curves $l$ meeting $F$ transversally we conclude 
$\varphi(F) \cong F$. The other assertions are evident.
\qed

\begin{proposition} \label{prop713}
Assume that $X$ has structure (O) and let $\varphi: X \to X'$ be a
birational Mori contraction. Then $\varphi$ contracts a divisor
$E \cong \bP_2$ with normal bundle $\sO_E(-1)$ to a point. $X'$ has
structure (O) and $\varphi_{|G}$ blows down one (--1)-curve in $G$.
\end{proposition}

\proof As $-K_F$ is divisible by two $\varphi$ is induced by a rational
curve $l$ with $-K_X \cdot l = 2$ and contracts $E \cong \bP_2$ to 
a point on the smooth threefold $X'$. Since $G \cdot l = 1$ and $l$ is
a line on the exceptional $\bP_2$ we have $G_{|E} = \sO_E(1)$. 
So $G^2 \cdot E = 1$ which also means that $E$ intersects $D = G^2$ (which
is an anticanonical divisor for $G$) in one point. Hence a general $l$ is 
a (--1)-curve in the appropriate $G$ containing it.
\qed

\subsubsection{The outcome -- Mori fiber spaces}

Since $\kappa(X)=-\infty$ the Mori program must terminate with a fiber space. 
We start with a smooth 3-fold $X$ having structure (A). After a finite number of
blow-downs described in \ref{prop710} and \ref{prop713} we obtain a smooth 3-fold 
(which we call again $X$) with structure (A) or (O). The final step in the Mori
program is a fiber type contraction $\varphi: X \to Y$.
We first consider the case where $\varphi$ is a conic bundle.

\begin{proposition}
Assume that $X$ has structure (A) and let $\varphi: X \to Y$ be a Mori
contraction to a smooth surface $Y$, i.~e.\ $\varphi$ is a conic bundle.
Let $l$ be a general conic. Then 
\begin{description}
\item either $F \cdot l = 1$ and $Y = F$, $X = \bP_1 \times F$ and $\varphi$
is the second projection. In this case $A = \varphi^* (C)$ for some elliptic
curve $C \in |-K_F|$.
\item or $F \cdot l = 0$. In this case $Y = \bP_1 \times \bP_1$ and $\varphi$
is a $\bP_1$-bundle; in particular $\varphi_{|F}$ gives $F$ the structure
of a ruled surface.
\end{description}
\end{proposition}

\proof a) We first consider the case $F \cdot l \neq 0$.
As $-K_X \cdot l = A \cdot l + 2 F \cdot l = 2$ the only possibility is $F \cdot l = 1$ 
which also implies that $\varphi$ is a $\bP_1$-bundle and $\varphi_{|F}$ 
is an isomorphism. Consider the product map $p = f \times \varphi: 
X \to \bP_1 \times Y$ which is generically one to one. If $D \subset X$ is a 
curve which is contracted by $p$ then $D$ is also contracted by $\varphi$. 
Therefore $D$ is a fiber of $\varphi$ and the rigidity lemma shows that 
$D$ cannot be contracted by $p$. It follows that $p$ is an isomorphism 
$X \cong \bP_1 \times F$ and $\varphi$ is the projection to the second 
factor.\\
b) The other case is $F \cdot l = 0$. This implies $F = \varphi^* F'$ 
which gives a factorization 
$f: X \stackrel{\varphi}{\la} Y \stackrel{pr_2}{\la} \bP_1$.
We notice that $pr_2 \circ \varphi_{|A}$ equals the projection of $A$ to 
the second factor. We also have $A \cdot l = 2$ which
means that $\varphi$ restricted to $A$ is generically 2:1. It is also
finite as $A = C \times \bP_1$ has no contractible curves. In fact if we look 
what happens for a general fiber $F$ then $\varphi_{|F}$ maps $C$ 2:1
on 
$pr_2^{-1}(f(F)) \cong \bP_1$. As $pr_2 \circ \varphi_{|A}$ is the second 
projection of $A$ the ramification locus of $\varphi_{|A}$ must be equal to 
some fibers of the first projection of $A$ which gives $Y = \bP_1 \times 
\bP_1$.\\
Let $Q = \{pt\} \times \bP_1$ be a general fiber of the first projection and 
let $X_Q = \varphi^{-1}(Q)$. Then $A_{|Q}$ consists of two sections 
$Q_1, Q_2$ of $\varphi_{|X_Q}$. As $A \cdot l = 2$ we conclude that 
$Q_i \cdot l = 1$ and $\varphi_{|X_Q}$ is a smooth conic bundle
i.~e.\ $X_Q$ is a ruled surface.(In fact $X_Q$ must be $\bP_1 \times \bP_1$ 
as $Q_1$ and $Q_2$ do not meet.) In particular the discriminant locus 
$\Delta$ of $\varphi$ is contained in some fibers of the first projection of 
$Y$. As $\varphi$ is an extremal contraction every nonsingular rational curve 
in $\Delta$ must meet the rest of $\Delta$ in at least two points (see for 
example \cite[p. 83]{Mi83}). In our situation this implies that $\Delta$ 
is empty and $\varphi$ is a $\bP_1$-bundle so $\varphi_{|F}$ exhibits 
$F = f^{-1}(pt)$ as a ruled surface over $pr_2^{-1}(pt)$.
\qed

\begin{proposition}
Assume that $X$ has structure (O) and let $\varphi: X \to Y$ be a Mori
contraction to a smooth surface $Y$. Then $Y \cong G$ and
$\varphi$ is a $\bP_1$-bundle.
\end{proposition}

\proof Let $l$ be a general conic. As $-K_X \cdot l = 2$ we get $G \cdot l = 1$
so $\varphi$ is regular and $G$ is a section.
\qed

Next we consider del Pezzo fibrations over some curve which must be $\bP_1$
as $X$ is rationally connected. Since $\varphi$ is a Mori contraction we also know $\rho(X)=2$.

\begin{proposition}
Assume that $X$ has structure (A) and let $\varphi: X \to \bP_1$ be a Mori
contraction. Then $\varphi = f$ and in particular $F$ is a del Pezzo surface.
\end{proposition}

\proof
The contraction is generated by a rational curve $l$ with $-K_X \cdot l \in 
\{1, 2, 3\}$ and $l$ is numerically effective. Therefore $A \cdot l \geq 0$
and in fact $A \cdot l > 0$ as $A$ is not a fiber of $\varphi$. If
$-K_X \cdot l = 1$ or $2$ this implies that $F \cdot l = 0$ so
$F = \varphi^*(pt)$ and $f$ and $\varphi$ coincide. In the last case
$-K_X \cdot l = 3$ and $\varphi$ is a $\bP_2$-bundle and $l$ a line. As
$\bP_2$ is not fibered, $F$ restricted to a $\bP_2$ is trivial and again
the two fibrations coincide.
\qed

\begin{proposition}
Assume that $X$ has structure (O) and let $\varphi: X \to \bP_1$ be a Mori
contraction. Then $\varphi$ is a quadric bundle and $\varphi_{|G}$ defines
a $\bP_1$-fibration on $G$.
\end{proposition}

\proof
As $-K_X$ is divisible by 2, $\varphi$ is always a quadric bundle. Let
$Q \cong \bP_1 \times \bP_1$ be a fiber. The canonical bundle formula
shows that $G_{|Q}$ has type $(1,1)$ so $Q$ intersects $G$ in a smooth
rational curve.
\qed

The last possibility is the case where the base of the Mori fiber space
is a point, i.~e.\ $X$ is a Fano threefold with $\rho = 1$.

\begin{proposition} 
Let $X$ be a Fano manifold with $\rho(X) = 1$. Then
$X$ does not have
structure (A).
\end{proposition}

\proof This is obvious because $F^2 = 0$ but $A_{|F}$ is the canonical
bundle of $F$ which is not zero.
\qed

\begin{proposition} Assume that $X$ is a Fano manifold with $\rho(X)=1$
and that $X$ has structure (O). Then either $X \cong \bP_3$ and $G$ is
a smooth quadric or $X$ has index two, $G$ is the generator of $\mathrm{Pic}
(X)$ and is therefore a del Pezzo surface of degree $1 \leq G^3 \leq 5$.
\end{proposition}

\proof We just use the fact that $-K_X$ is divisible by two and cite
Iskovskih's classification.
\qed

\vspace{1cm} 
\small 
\begin{tabular}{lcl} 
Thomas Bauer and Thomas Peternell\\ Mathematisches Institut \\ Universit\" at Bayreuth \\ D-95440 Bayreuth, Germany \\ thomas.bauer@uni-bayreuth.de \\ 
thomas.peternell@uni-bayreuth.de  \\ 
\end{tabular}


\begin{thebibliography}{8authors}

\bibitem[8authors]{8authors} T.~Bauer, F.~Campana, T.~Eckl, S.~Kebekus, 
T.~Peternell, 
S.~Rams, T.~Szemberg, L.~Wotzlaw: \emph{A reduction map for nef line bundles.}
Complex Geometry (Collection of papers dedicated to Hans Grauert), 
Springer-Verlag, Berlin Heidelberg New York 2002, 27--36

\bibitem[BPV84]{BPV84}W.~Barth, C.~Peters, A.~van de Ven: \emph{Compact complex 
surfaces.} Ergebnisse der Mathematik und ihrer Grenzgebiete (3), 
Springer-Verlag, Berlin 1984

\bibitem[BS95]{BS95}M.~Beltrametti, A.J.~Sommese:  \emph{The adjunction theory of complex projective varieties.} de Gruyter 1995 

\bibitem[DPS93]{DPS93} J.-P.~Demailly, T.~Peternell, M.~Schneider:
\emph{K\"ahler manifolds with numerically effective Ricci class.}
Compositio Math.\ 89 (1993), 217--240

\bibitem[DPS94]{DPS94} J.-P.~Demailly, T.~Peternell, M.~Schneider:
\emph{Compact complex manifolds with numerically effective tangent bundles.}
J.\ Algebraic Geom.\ 3 (1994), 295--345

\bibitem[DPS96]{DPS96} J.-P.~Demailly, T.~Peternell, M.~Schneider:
\emph{Compact K\"ahler manifolds with Hermitian semipositive anticanonical
bundle.} Compositio Math.\ 101 (1996), 217--224

\bibitem[DPS01]{DPS01} J.-P.~Demailly, T.~Peternell, M.~Schneider: 
\emph{Pseudo-effective line bundles on compact K\"ahler manifolds.}
Internat.\ J.\ Math.\ 12 (2001), 689--741

\bibitem[Fu90]{Fu90} T.~Fujita: \emph{Classification theories of polarized varieties.}  London Math. Soc. Lect. Notes Ser. 155 (1990) 

\bibitem[Ka85]{Ka85} Y.~Kawamata: \emph{Pluricanonical systems on
minimal algebraic varieties.} Invent.\ Math.\ 79 (1985), 567--588

\bibitem[KM98]{KM98} J.~Koll\'ar, S.~Mori: \emph{Birational geometry of algebraic 
varieties.} Cambridge Tracts in Mathematics 134, Cambridge University Press, 
Cambridge 1998


\bibitem[Mo82a]{Mo82a} S.~Mori: \emph{Threefolds whose canonical bundles are 
not
numerically effective.} Ann.\ of Math.\ (2) 116 (1982), 133--176

\bibitem[Mo82b]{Mo82b} S.~Mori: \emph{Threefolds whose canonical bundles are 
not numerically effective.} Algebraic threefolds (Varenna 1981),
Lecture Notes in Math.\ 947, Springer-Verlag, Berlin-New York 1982, 155--189

\bibitem[Mi83]{Mi83} M.~Miyanishi: \emph{Algebraic methods in the theory of
algebraic threefolds.} Algebraic varieties and analytic varieties
(Tokyo 1981), Adv.\ Stud.\ Pure Math.\ 1, North-Holland, Amsterdam 1983, 69--99

\bibitem[Na87]{Na88} N.~Nakayama: \emph{On Weierstrass models. Alg.Geom. and Comm. Algebra}  vol. in
honour of Nagata, vol. 2. Kinokuniya, Tokyo 1988, 405--431

\bibitem[PS98]{PS98} T.~Peternell, F.~Serrano: \emph{Threefolds with nef
anticanonical bundles.} Collect.\ Math.\ 49 (1998), 465--517

\bibitem[Sa82]{Sa82} F.~Sakai: \emph{Anti-Kodaira dimension of ruled surfaces.}
Sci.~Rep., Saitama Univ. X,2 (1982), 1--7


\bibitem[Zh96]{Zh96} Q.~Zhang: \emph{On projective manifolds with nef anticanonical bundle.} 
J.f.d.\ reine u.\ angew.\ Math.\ 478 (1996) 57--60

\end{thebibliography}
\end{document}